\documentclass[10pt]{amsart}

\newtheorem{theorem}{Theorem}
\newtheorem{proposition}{Proposition}
\newtheorem{lemma}{Lemma}
\newtheorem{corollary}{Corollary}
\newtheorem{assumption}{Assumption}

\newtheorem{theoremA}{Watanabe's Theorem}

\theoremstyle{definition}
\newtheorem{definition}{Definition}

\theoremstyle{remark}
\newtheorem{remark}{Remark}

\newenvironment{pf}{\proof}{\endproof \medskip}
\newenvironment{pf*}[1]{\medskip \noindent {\em #1.} }{\endproof \medskip}

\newcommand{\xbrace}[1]{\langle #1 \rangle}

\newcommand{\xnorm}[1]{ \Vert #1 \Vert }

\newcommand{\zz}[1]{\mathbb #1}

\title{Spatial Epidemics: Critical Behavior in One Dimension} 
\author{Steven P. Lalley} 

\thanks{Supported by NSF grant DMS-0405102}
\date{\today}

\begin{document}

\begin{abstract}
In the simple mean-field \emph{SIS} and \emph{SIR} epidemic models,
infection is transmitted from infectious to susceptible members of a
finite population by independent $p-$coin tosses. Spatial variants of
these models are considered, in which finite populations of size $N$ are
situated at the sites of a lattice and infectious contacts are limited
to individuals at neighboring sites. Scaling laws for these models are
given when the infection parameter $p$ is such that the epidemics are
\emph{critical}. It is shown that in all cases there is a
\emph{critical threshold} for the numbers initially infected: below
the threshold, the epidemic evolves in essentially the same manner as
its \emph{branching envelope}, but at the threshold evolves like a
branching process with a size-dependent drift. The corresponding
scaling limits are super-Brownian motions and Dawson-Watanabe
processes with killing, respectively.
\end{abstract}

\maketitle 

\section{Introduction}
\label{sec:intro}

\subsection{Critical mean-field epidemics: threshold
behavior}\label{ssec:criticalMF} It was discovered by Martin-L\"{o}f
\cite{martin-lof} and independently by Aldous \cite{aldous} that the
simple mean-field \emph{SIR} epidemic, also known as the
\emph{Reed-Frost} epidemic, exhibits a curious threshold behavior at
criticality. Roughly, if $U_{N}$ is the \emph{size} -- that is, the
number of individuals ever infected -- of the epidemic in a population
of size $N$, then $U_{N}$ has a markedly different asymptotic
distribution when the number $J_{0}$ of individuals infected at time
$0$ is of order $N^{1/3}$ than when it is of order $o (N^{1/3})$. In
particular, if $J_{0}\sim bN^{\alpha}$ as $N \rightarrow \infty$, then
\begin{equation}\label{eq:sizeAsymptotics}
	U_{N}/N^{2\alpha}  \stackrel{\mathcal{D}}{\longrightarrow} 
	\tau_{b} ,
\end{equation}
where $\tau_{b}$ is the first passage time to the level $b$ by a
standard Wiener Process, if $\alpha <1/3$, or by a Wiener process with
time-dependent drift $t$, if $\alpha =1/3$.  This reflects the fact
that the size of the largest connected component in a critical
($p=1/N$) Erd\"{o}s-Renyi random graph on $N$ vertices is of order
$N^{1/3}$. R.~Dolgoarshinnykh and the author
\cite{dolgoarshinnykh-lalley} have observed that there is a similar
critical threshold effect for the simple mean-field \emph{SIS}
epidemic, but here the threshold for $J_{0}$ is at $N^{1/2}$ rather
than $N^{1/3}$, and the limit distribution at the threshold involves
first passage times by Ornstein-Uhlenbeck processes. In fact, there is
an asymptotic form for the entire evolution of the epidemic at
criticality that undergoes a discontinuity at $J_{0}\approx N^{1/2}$:
If $J_{n}$ denotes the number of individuals infected at time $n$ then
\begin{equation}\label{eq:SIS-mf}
	N^{-\alpha }J_{[N^{\alpha}t]}  \stackrel{\mathcal{D}}{\longrightarrow} Y_{t}
\end{equation}
where $Y_{0}=b$ and $Y_{t}$ is either a Feller diffusion or a Feller
diffusion with location-dependent drift $-Y_{t}^{2} \, dt$, depending
on whether $\alpha <1/2$ or $\alpha = 1/2$, that is,
\begin{align}\label{eq:sisCriticalScaling}
	dY_{t} &=\sqrt{Y_{t}} dW_{t} \quad \text{if} \; \alpha <1/2;\\
\notag	
	 dY_{t} &=-Y_{t}^{2}\,dt +\sqrt{Y_{t}} dW_{t} \quad \text{if} \; \alpha =/2.
\end{align}
Note that in the case $\alpha <1/2$ the limit process -- the Feller
diffusion -- is the same as the limit process for the rescaled
critical Galton-Watson process. There is a similar process-level
threshold effect for the Reed-Frost epidemic at $J_{0}\approx N^{1/3}$
-- see \cite{dolgoarshinnykh-lalley}.

\subsection{Spatial \emph{SIS} and \emph{SIR}
epidemics}\label{ssec:spatialSIRS} The purpose of this article is to
show that there is a similar critical threshold effect for spatial
epidemics in one spatial dimension -- see Theorem~\ref{theorem:sis}
below. The epidemic models considered are  simple discrete-time
spatial  analogues of the Reed-Frost and stochastic logistic
epidemics. These are chosen primarily to streamline the mathematical
analysis; however,  analogous effects should
also be expected for more complex models in dimension $d=1$. A
secondary motivation for the specification of the  spatial \emph{SIR}
model is that it has a percolation (random graph) description similar to the
Erd\"{o}s-Renyi random graph description of the Reed-Frost epidemic,
and so our main result can be interpreted as a statement about the
connected clusters in certain percolation models.  

The spatial epidemics, which we will henceforth call the
\emph{SIR}$-d$ and \emph{SIS}$-d$ epidemics, are defined as follows:
Assume that at each lattice point $x\in \zz{Z}^{d}$ is a homogeneous
population (``village'') of $N$ individuals, each of whom may at any
time be either susceptible or infected, or (in the \emph{SIR}
variants) recovered. As in the corresponding mean-field models (see
\cite{martin-lof}), infected individuals remain infected for one unit
of time, and then recover; in the \emph{SIR}-$d$ epidemic, infected
individuals recover and are thereafter immune from infection, while in
the \emph{SIS-}$d$ model, infected individuals, upon recovery, become
once again susceptible to infection. The rules of infection are the
same as for the corresponding mean-field models, \emph{except} that
the infection rates depend on the locations of the infected and
susceptible individuals. Thus, at each time $t=0,1,2,\dotsc $, for
each pair $(i_{x},s_{y})$ of an infected individual located at $x$ and
a susceptible individual at $y$, the disease spreads from $i_{x}$ to
$s_{y}$ with probability $p_{N} (x,y)$. We shall only consider the
case where the transmission probabilities $p_{N} (x,y)$ are spatially
homogeneous, nearest-neighbor, and symmetric, and scale with the
village size $N$ in such a way that the expected number of infections
by a contagious individual in an otherwise healthy population is 1 (so
that the epidemic is \emph{critical}), that is,
\begin{assumption}\label{assumption:PoissonScaling}
\begin{align}\label{eq:tps}
	p_{N} (x,x+e_{i})	&= C_{d}/N \quad \text{if}
	\quad |e_{i}|=1 \; \text{or} \; 0,  \quad \text{where} \\
	\quad C_{d}&=1/ (2d+1).
\end{align}
\end{assumption}

Similar models incorporating separated clusters have been studied by
Schinazi \cite{schinazi}, Belhadji \& Lanchier
\cite{belhadji-lanchier}, and others, but these studies have focused
on \emph{SIS} variants of the models where all infected individuals in
a colony recover simultaneously, and where infection rates within and
between colonies vary.  Critical behavior of certain spatial epidemic
models has been addressed in the literature, in particular for
long-range contact processes \cite{mueller-tribe},
\cite{durrett-perkins}, which are in certain respects similar to the
\emph{SIS-d} model described above; however, critical behavior of
spatial \emph{SIR} models has not been previously studied. For surveys
of contact models in spatial epidemics, see \cite{mollison} and
\cite{durrett:survey}. 

Interest in spatial epidemic models has largely been focused on
dimensions $d\geq 2$, and especially $d=2$, for natural
reasons. Nevertheless, nearest neighbor infection models in dimension
$d=1$ may be of interest in certain contexts: Many plant and animal
species live in river valleys or along shorelines, and for these the
natural  dimension  for spatial interactions is $d=1$.

\subsection{Epidemic Models and Random Graphs}\label{ssec:rgraphs} The
models described above have equivalent descriptions as structured
random graphs, that is, percolation processes. Consider first the
simple \emph{SIR} (Reed-Frost) epidemic. In this model, no individual
may be infected more than once; furthermore, for any pair $x,y$ of
individuals, there will be at most one opportunity for infection to
pass from $x$ to $y$ or from $y$ to $x$ during the course of the
epidemic. Thus, one could simulate the epidemic by first tossing a
$p-$coin for every pair $x,y$, drawing an edge between $x$ and $y$ for
each coin toss resulting in a Head, and then using the resulting
(Erd\"{o}s-Renyi) random graph determined by these edges to determine
the course of infection in the epidemic. In detail: If $Y_{0}$ is the
set of infected individuals at time $0$, then the set $Y_{1}$ of
individuals infected at time $1$ consists of all $x\not \in Y_{0}$
that are connected to individuals in $Y_{0}$, and for any subsequent
time $n$, the set $Y_{n+1}$ of individuals infected at time $n+1$
consists of all $x\not \in \cup_{j\leq n} Y_{j}$ who are connected to
individuals in $Y_{n}$. Note that the set of individuals ultimately
infected during the course of the epidemic is the union of those
connected components of the random graph containing at least one
vertex in $Y_{0}$.

Similar random graph descriptions may be given for the mean-field
\emph{SIS} and the spatial \emph{SIS} and \emph{SIR} epidemic
models. Consider for definiteness the \emph{SIR-d} epidemic. To
simulate this, first build a random graph by Bernoulli bond
percolation on the graph $\zz{K}_{N}\times \zz{Z}^{d}$, where
$\zz{K}_{N}$ is the complete graph on $N$ vertices. Given the random
graph, simulate the generations $Y_{n}$ of the \emph{SIR-d} epidemic
by the same rule as in the mean-field case: For each generation $n$,
define the set $Y_{n+1}$ of individuals infected at time $n+1$ to be
the set of all vertices $x\not \in \cup_{j\leq n} Y_{j}$ who are
connected to individuals in $Y_{n}$. Similar random graph descriptions
may be given for \emph{SIS} epidemics, but using oriented percolation
for the  random graphs.

\subsection{Branching envelopes of spatial epidemics}
\label{ssec:breSpatial}

The \emph{branching envelope} of a spatial \emph{SIS}$-d$ or
\emph{SIR}$-d$ epidemic is a \emph{branching random walk} on the
integer lattice $\zz{Z}^{d}$. This evolves as follows: Any particle
located at site $x$ at time $t$ lives for one unit of time and then
reproduces, placing random numbers $\xi_{y}$ of offspring at the sites
$y$ such that $|y-x|\leq 1$. The random variables $\xi_{y}$ are
mutually independent, each with Binomial-$(N,C_{d}/N)$ distributions,
where $C_{d}=1/ (2d+1)$. Denote this reproduction rule by
$\mathcal{R}_{N}$, and denote by $\mathcal{R}_{\infty}$ the
corresponding offspring law in which the Binomial distribution is
replaced by the Poisson distribution with mean $C_{d}$.  Note that for
each of the offspring distributions $\mathcal{R}_{N}$, the branching
random walk is critical, that is, the expected total number of
offspring of a particle is $1$.

A fundamental theorem of S.~Watanabe \cite{watanabe} asserts that,
under suitable rescaling (the \emph{Feller scaling}) the
measure-valued processes naturally associated with critical branching
random walks converge to a limit, the \emph{standard Dawson-Watanabe}
process, also known as \emph{super Brownian motion}. 


\begin{definition}\label{definition:FellerScaling}
The \emph{Feller-Watanabe scaling operator} $\mathcal{F}_{k}$ scales
mass by $1/k$ and space by $1/\sqrt{k}$, that is, for any finite Borel
measure $\mu (dx)$ on $\zz{R}^{d}$ and any test function $\phi (x)$,
\begin{equation}\label{eq:FellerScaling}
	\xbrace{\phi ,\mathcal{F}_{k} \mu }=
	k^{-1}\int \phi (\sqrt{k}x) \mu (dx)
\end{equation}
\end{definition}

\begin{theoremA}\label{theorem:watanabe}
Fix $N$, and for each $k=1,2,\dotsc$ let $Y^{k}_{t}$ be a branching
random walk with offspring distribution $\mathcal{R}_{N}$ and initial
particle configuration $Y^{k}_{0}$. (In particular, $Y^{k}_{t} (x)$
denotes the number of particles at site $x\in \zz{Z}$ in generation
$[t]$, and $Y^{k}_{t}$ is the corresponding Borel measure on $\zz{R}$.)
If the initial mass distributions converge, after rescaling, as $k
\rightarrow \infty$, that is, if
\begin{equation}\label{eq:initialCondition}
	\mathcal{F}_{k}Y^{k}_{0}
	\Longrightarrow X_{0}
\end{equation}
for some finite Borel measure $X_{0}$ on $\zz{R}^{d}$, then the
rescaled measure-valued processes $\mathcal{F}_{k}Y^{k}_{kt}$ converge
in law as $k \rightarrow \infty$:
\begin{equation}\label{eq:watanabeConclusion}
	(\mathcal{F}_{k}Y^{k})_{kt} \Longrightarrow X_{t}.
\end{equation}
\end{theoremA}

The limit is  the \emph{standard Dawson-Watanabe process} $X_{t}$
(also known as \emph{super-Brownian motion}). See \cite{etheridge} for
more.  In dimension $d=1$ the random measure
$X_{t}$ is for each $t$ absolutely continuous relative to Lebesgue
measure \cite{konno-shiga}, and the Radon-Nikodym derivative $X(t,x)$
is jointly continuous in $t,x$ (for $t>0$). In dimensions $d\geq 2$
the measure $X_{t}$ is almost surely singular, and is supported by a
Borel set of Hausdorff dimension $2$ \cite{dawson-hochberg}. 

\subsection{Scaling limits of critical spatial epidemics}
\label{ssec:scalingLimits} Our main result, Theorem \ref{theorem:sis}
below, asserts that after appropriate rescaling the \emph{SIS} -$1$
and \emph{SIR} -$1$ spatial epidemics converge weakly as $N
\rightarrow \infty$.  The limit processes are either standard
Dawson-Watanabe processes or {Dawson-Watanabe processes with
variable-rate killing}, depending on the  initial configuration of
infected individuals.  The Dawson-Watanabe process $X_{t}$ with
killing rate $\theta=\theta (x,t,\omega )$ (assumed to be
progressively measurable and jointly continuous in $t,x$) and variance
parameter $\sigma^{2}$ can be characterized by a martingale problem
(\cite{dawson-perkins}, sec. 6.2):
for every $\phi \in C^{2}_{c} (\zz{R}^{d})$,
\begin{equation}\label{eq:dwMG}
	\xbrace{X_{t},\phi}-\xbrace{X_{0},\phi}-
	\frac{\sigma_{d}}{2}\int_{0}^{t} \xbrace{X_{s},\Delta \phi}
	\,ds +
	\int_{0}^{t} \xbrace{X_{s},\theta (\cdot ,s)\phi} \,ds
\end{equation}
is a martingale with the same quadratic variation as for the standard
Dawson-Watanabe process. Existence and distributional uniqueness of such
processes in general is asserted in \cite{dawson-perkins} and proved,
in various cases, in \cite{dawson:geo} and \cite{evans-perkins}. It is
also proved in these articles that the law of a Dawson-Watanabe
process with killing on a finite time interval is absolutely
continuous with respect to that of a standard Dawson-Watanabe process
with the same variance parameter, and that the likelihood ratio
(Radon-Nikodym derivative) is \cite{dawson-perkins}
\begin{equation}\label{eq:dwLR}
	\exp \left\{-\int \theta (t,x) \,dM (t,x)-\frac{1}{2} 
	\int \xbrace{X_{t},\theta (t,\cdot)^{2}}\, dt \right\} ,
\end{equation}
where $dM (t,x)$ is the \emph{orthogonal martingale measure} attached
to the standard Dawson-Watanabe process (see \cite{walsh} and
sec.~\ref{ssec:omm} below).  Absolute continuity implies that sample
path properties are inherited: In particular, if $X_{t}$ is a
one-dimensional Dawson-Watanabe process with killing, then almost
surely the random measure $X_{t}$ is absolutely continuous, with
density $X (x,t)$ jointly continuous in $x$ and $t$.

\begin{theorem}\label{theorem:sis}
Let $Y^{N}_{t} (x)$ be the number infected at time $t$ and site $x$ in
a critical \emph{SIS}-$1$ or  \emph{SIR}-$1$ epidemic with village
size $N$ and initial configuration $Y^{N}_{0} (x)$.
 Fix $\alpha >0$, and let $X^{N} (t,x)$ be the renormalized particle
 density function  process obtained  by linear interpolation in $x$ from the
 values 
\begin{equation}\label{eq:epiDensityProcess}
	X^{N} (t,x):=\frac{Y^{N}_{[N^{\alpha}t]} (\sqrt{N^{\alpha }}
	x)}{\sqrt{N^{\alpha}}}
	\quad \text{for} \quad x\in \zz{Z}/\sqrt{N^{\alpha}} .
\end{equation}
Assume that there is a compact interval $J$ such that the initial
particle density functions $X^{N} (0,x)$ all have support contained in $J$,
and assume that the functions $X^{N} (0,x)$ converge in $C_{b} (\zz{R})$
to a function $X (0,x)$. Then under Assumption
\ref{assumption:PoissonScaling}, as $N \rightarrow \infty$,
\begin{equation}\label{eq:limitSIS}
	X^{N} (t,x) \Longrightarrow X (t,x)
\end{equation}
where $X (t,x)$ is the density of a Dawson-Watanabe process $X_{t}$
with initial density $X (0,x)$ and killing rate $\theta $ depending on the value
of $\alpha$ and the type of epidemic (SIS or SIR) as follows:
\begin{enumerate}
\item [(a)] \emph{SIS:} If $\alpha <\frac{2}{3}$ then $\theta (x,t)=0$.
\item [(b)] \emph{SIS:}  If $\alpha =\frac{2}{3}$ then
\begin{equation}\label{eq:sisKillingRate}
	\theta (x,t)=X (x,t)/2.
\end{equation}
\item [(c)]   \emph{SIR:}  If $\alpha <\frac{2}{5}$ then then $\theta (x,t)=0$.
\item [(d)]  \emph{SIR:}  If $\alpha =\frac{2}{5}$ then
\begin{equation}\label{eq:sirKillingRate}
	\theta (x,t)=\int_{0}^{t} X (x,s)\, ds.
\end{equation}
\end{enumerate}
The convergence $\Rightarrow$ in \eqref{eq:limitSIS} is weak
convergence relative to the Skorohod topology on the space
$\zz{D}([\cup ,\infty ),C_{b} (\zz{R}))$ of  cadlag functions
$t\mapsto X(t,x)$  valued in $C_{b} (\zz{R})$.
\end{theorem}

The proof of Theorem \ref{theorem:sis} is given in
section~\ref{sec:spatial}; it will will depend on
Theorem~\ref{theorem:BRW} below. In both Theorems \ref{theorem:BRW}
and \ref{theorem:sis}, the assumption that the initial particle
densities have common compact support $J$ can undoubtedly be weakened;
however, this assumption eliminates certain technical complications in
the arguments (see equation \eqref{eq:maxSpatialBRW} in
sec.~\ref{sec:spatial}).

\medskip \noindent \textbf{Remarks.}  (A) The case of principal
interest is the \emph{SIR-1}  epidemic. The \emph{SIS-1} epidemic is
closely related to the \emph{long-range
contact process} studied by Mueller and Tribe \cite{mueller-tribe} in
$d=1$ and by Durrett and Perkins \cite{durrett-perkins} in $d\geq 2$,
and in particular the limit process for the \emph{SIS-1} epidemic, the
Dawson-Watanabe process with killing rate \eqref{eq:sisKillingRate},
is the same as that for the rescaled contact process in
$d=1$. Nevertheless, the rescaled contact process and the
discrete-time \emph{SIS-1} epidemic studied here differ in some
important technical respects, so parts (a)- (b) of
Theorem~\ref{theorem:sis} do not follow from the results of
\cite{mueller-tribe}. 

\medskip (B) The proof strategy here (sec.~\ref{sec:spatial} below) is
considerably different -- and simpler -- than the usual approach,
based on martingale methods, taken in the literature of weak
convergence to superprocesses, such as in \cite{mueller-tribe},
\cite{durrett-perkins}, \cite{cox-durrett-perkins}, and
\cite{durrett-mytnik-perkins}. Although martingale methods might be
made to work here, the prospect of using them in connection with
processes such as the \emph{SIR-1} epidemic with history-dependent
rates is rather daunting. Instead, we rely on the fact that the laws
of both the \emph{SIS} and \emph{SIR} epidemics are absolutely
continuous with respect to those of their branching envelopes, for
which scaling limits are already known. The Radon-Nikodym derivatives
have tractable forms, as exponentials of certain stochastic
integrals. These will be shown to converge to the corresponding
Radon-Nikodym derivatives \eqref{eq:dwLR}. The advantage of this
strategy is that, given Theorem~\ref{theorem:BRW}, there is no need to
check tightness for the rescaled epidemic processes. Moreover, there
is virtually no additional work involved in establishing the result
for the \emph{SIS} model -- all that is needed is an additional simple
asymptotic estimation of the Radon-Nikodym derivatives.

\medskip 
(C) In higher dimensions, there is no analogous threshold effect at
criticality. In dimensions $d\geq 3$, a branching random walk started
by any number of particles will quickly diffuse, so that after a short time most 
occupied sites will have only $O (1)$ particles. Consequently, in both the
\emph{SIS-}$d$ and the \emph{SIR-} $d$ epidemics, the effect
of finite population size on the production of new infections will
be limited, and so both epidemics will behave in more or less the same
manner as their branching envelopes. In dimension $d=2$, the situation
is somewhat more interesting: it appears that here finite population
size will manifest itself by a logarithmic drag on the production of
new particles. The cases $d\geq 2$ will be discussed in detail in a
forthcoming paper of Xinghua Zheng \cite{zheng:thesis}.

\subsection{Heuristics: The standard coupling}\label{ssec:heuristics}
The critical thresholds for the \emph{SIS}$-d$ and \emph{SIR}$-d$
epidemics can be guessed by a simple comparison argument based on the
\emph{standard coupling} of the epidemic and its branching
envelope. For the \emph{SIS}$-d$ epidemic, the coupling is constructed
as follows: Build a branching random walk whose initial state
coincides with that of the epidemic, with particles to be colored
\emph{red} or \emph{blue} according to whether or not they represent
infections that actually take place (red particles represent actual
infections). Initially, all particles are red. At each time
$t=0,1,2,\dotsc$, particles of the branching random walk produce
offspring at neighboring sites according to the law described in
sec.~\ref{ssec:breSpatial} above. Offspring of blue particles are
always blue, but offspring of red particles may be either red or blue,
with the choices made as follows: All offspring of red particles at a
location $y$ choose labels $j\in [N]$ at random, independently of all
other particles; for any label $j$ chosen by $k>1$ particles, one of the
particles is chosen at random and colored red, and the remaining $k-1$
particles are colored blue. The population of all particles evolves as
a branching random walk, by construction, while the subpopulation of
red particles evolves as an \emph{SIS}$-d$ epidemic. Observe that the
branching random walk dominates the epidemic: thus, the duration,
size, and spatial extent of the epidemic are limited by those of the
branching envelope.

The standard coupling of an \emph{SIR}$-d$ epidemic with its branching
envelope is constructed in a similar fashion, but with the following
rule governing choices of color by offspring of red particles: All
offspring of red particles at a location $y$ choose numbers $j\in [N]$
at random, independently of all other particles.  If a particle
chooses a number $j$ that was previously chosen by a particle of an
earlier generation at the same site $y$, then it is assigned color
blue. If $k>1$ offspring of red particles choose the same number $j$
at the same time, and if $j$ was not chosen in an earlier generation,
then $1$ of the particles is assigned color red, while the  remaining
$k-1$ are assigned color blue. Under this rule, the subpopulation of
red particles evolves as an \emph{SIR}$-d$ epidemic.

In both couplings, the production of blue offspring by red particles
may be viewed as an attrition of the red population. Assume that
initially there are $N^{\alpha}$ particles; then by Feller's limit
theorem for critical Galton-Watson processes, the branching envelope
can be expected to survive for $O_{P} (N^{\alpha})$ generations, and
at any time prior to extinction the population will have $O_{P}
(N^{\alpha})$ members. These will be distributed among the sites at
distance $O_{P} (N^{\alpha /2})$ from the origin, and therefore in
dimension $d=1$ there should be about $O_{P} (N^{\alpha /2})$
particles per site. Consequently, for the \emph{SIS}$-1$ epidemic, the
rate of attrition per site per generation should be $O_{P}
(N^{\alpha-1})$, and so the total attrition rate per generation should
be $O_{P} (N^{3\alpha /2 -1})$. If $\alpha =2/3$, then the total
attrition rate per generation will be $O_{P} (1)$, just enough so that
the total attrition through the duration of the branching random walk
envelope will be on the same order of magnitude as the population size
$N^{\alpha}$.

For the \emph{SIR}$-1$ epidemic there is a similar heuristic
calculation. As for the \emph{SIS}$-1$ epidemic, the branching
envelope will survive for $O_{P} (N^{\alpha})$ generations, and up to
the time of extinction the population should have $O_{P} (N^{\alpha})$
individuals, about $O_{P} (N^{\alpha /2})$ per site. Therefore,
through $N^{\alpha}$ generations, about $N^{\alpha}\times N^{\alpha
/2}$ numbers $j$ will be retired, and so the attrition rate per site
per generation should be $O_{P} (N^{\alpha /2}\times N^{3\alpha /2})$,
making the total attrition rate per generation $O_{P} (N^{5\alpha
/2})$. Hence, if $\alpha =2/5$ then the total attrition per generation
should be $O_{P} (1)$, just enough so that the total attrition through
the duration of the branching random walk envelope will be on the same
order of magnitude as the population size.

\subsection{Weak convergence in $\zz{D} ([0,\infty ),C_{b}
(\zz{R}))$}\label{ssec:StrongWatanabe} The heuristic argument above
has an obvious gap: it relies crucially on the assertion that the
particles of a critical branching random walk distribute themselves
somewhat uniformly, at least locally, among the sites at distance
$N^{\alpha /2}$ from the origin.  The fact that the Dawson-Watanabe
process in dimension one has a \emph{continuous} density suggests that
this should be true, but does not imply it. Following is a
strengthening of the Watanabe theorem suitable for our purposes.

Denote by $C_{b} (\zz{R})$ the space of continuous, bounded,
real-valued functions on $\zz{R}$ with the sup-norm topology, and by
$\zz{D} ([0,\infty ),C_{b} (\zz{R}))$ the Skorohod space of cadlag
functions $X (t,x)$ valued in $C_{b} (\zz{R})$ (thus, for each $t\geq
0$ the function $X (t,x)$ is a continuous, bounded function of $x$).
Fix $N\leq \infty$, and for $k=1,2,\dotsc$ let $Y^{k}_{t} (x)$ be the
number of particles at site $x$ at time $[t]$ in a branching random
walk with offspring distribution $\mathcal{R}_{N}$ and initial
particle configuration $Y^{k}_{0} (\cdot)$. Let $X^{k} (t,x)$ be the
renormalized density function: that is, the function obtained by
linear interpolation (in $x$) from the values
\begin{equation}\label{eq:Yn}
	X^{k} (t,x)=\frac{Y^{k}_{[kt]} (\sqrt{k}x)}{\sqrt{k}}
	\quad \text{for} \quad x\in \zz{Z}/\sqrt{k}.
\end{equation}

\begin{theorem}\label{theorem:BRW}
Assume that there is a compact interval $J\subset \zz{R}$ such that
all of the the initial particle densities $X^{k} (0 ,\cdot)$ have
support $\subset J$, and assume that as $k \rightarrow \infty $ the
functions $X^{k} (0 ,\cdot)$ converge in $C_{b}(\zz{R})$ to a continuous
function $X (0,\cdot)$.  Then as $k \rightarrow \infty$,
\begin{equation}\label{eq:strongWatanabe}
	X^{k} (t,x) \Longrightarrow X (t,x),
\end{equation}
where $X (t,x)$ is the density function of a Dawson-Watanabe process
with initial density $X (0,x)$, and $\Rightarrow$ indicates weak
convergence relative to the Skorohod topology on $\zz{D} ([0,\infty
),C_{b} (\zz{R}))$. 
\end{theorem}

To prove this, it suffices to show that the sequence $X^{k} (t,x)$ of
densities is tight in $\zz{D} ([0,\infty ),C_{b} (\zz{R}))$, because
Watanabe's theorem implies that any weak limit of a subsequence must
be a density of the Dawson-Watanabe process. The proof of tightness,
carried out in section~\ref{sec:tightness} below, will be based on a
form of the Kolmogorov-Chentsov tightness criterion and a moment
calculation. This proof will use only three properties of the
offspring law $\mathcal{R}_{N}$:
\begin{itemize}
\item [(a)] the mean number of offspring is $1$;
\item [(b)]  the total number of offspring has
finite $m$th moment, for each $m\geq 1$; and 
\item [(c)]  offspring choose locations at random from among the
neighboring sites.
\end{itemize}
Since the $m$th moments are bounded uniformly
over the class of Binomial-$(N,1/N)$ distributions, it will follow that
tightness holds simultaneously for all of the offspring laws
$\mathcal{R}_{N}$. Only the case $N=\infty$ will be needed for the
analysis of the spatial epidemic models, however.

\medskip \noindent \textbf{Remark.} Mueller and Tribe
\cite{mueller-tribe} proved that the density processes of rescaled
long-range contact processes in one dimension converge weakly to the
density process of a Dawson-Watanabe process with killing, but make no
explicit use of branching random walks in their argument.
Nevertheless, they likely were aware that the density processes of
rescaled branching random walks would also converge weakly in one
dimension.


%

\subsection{Spatial extent of the Dawson-Watanabe
process}\label{ssec:weierstrass} 

An object of natural interest in connection with the spatial SIS and
SIR epidemics is the \emph{spatial extent} of the process, that is,
the area reached by the infection. Under certain natural restrictions
on the initial configuration of infected individuals, the spatial
extent will, by Theorem \ref{theorem:sis}, be well-approximated in law
by the area covered by the limiting Dawson-Watanabe process, after
suitable scaling. For Dawson-Watanabe processes $X_{t}$ with
location-dependent killing, the distribution of $\mathcal{R}
(X):=\cup_{t}\text{support} (X_{t})$ likely cannot be described in
closed form. However, for the Dawson-Watanabe process with constant
killing rate, the distribution of $\mathcal{R} (X)$ can be given in a
computable form, as we now show.

\begin{proposition}\label{proposition:weierstrass}
Let $X_{t}$ be the standard one-dimensional Dawson-Watanabe process
with variance parameter $\sigma^{2}=1$.  For any finite Borel measure
$\mu$ with support contained in the interval $D= (0,a)$,
\begin{equation}\label{eq:weierstrass}
		-\log P (\mathcal{R} (X)\subset D \, | \, X_{0}=\mu )
		=\int\wp_{L} (x/\sqrt{6})\,\mu (dx)
\end{equation}
where  
\begin{equation}\label{eq:P}
	\wp_{L} (x)=\sum_{\omega \in L\setminus\{0\}}
			\left\{\frac{1}{6
			(x-\omega)^{2}}-\frac{1}{\omega^{2}} 
			\right\}
\end{equation}
is the Weierstrass $\wp -$function with period lattice $L$ generated
by $\sqrt{6}ae^{\pi i/3}$.
\end{proposition}

\begin{pf}
By a theorem of Dynkin \cite{etheridge}, ch.~8, the function 
\begin{equation}\label{eq:ud}
	u_{D} (x):=-\log  P (\mathcal{R} (X)\subset D \, | \, X_{0}=\delta_{x} )
\end{equation}
is the unique solution of the differential equation $u''=u^{2}$ in $D$
with boundary conditions $u (x)\rightarrow \infty$ as $x \rightarrow
0,a$. Set $v=u'$; then the  equation $u''=u^{2}$ becomes $v'=u^{2}$,
and so $dv/du =u^{2}/v$. Integration gives
\[
	v^{2}/2=u^{3}/3 +C,
\]
which, up to constants, is the differential equation of the
$\wp-$function \cite{whittaker-watson}, ch.~20. The result
\eqref{eq:weierstrass} now follows for the special case $\mu
=\delta_{x}$. The general case now follows by the superposition
principle for branching processes: in particular, the Dawson-Watanabe
process $X_{t}$ with initial condition $X_{0}=\mu +\nu $ can be
decomposed as the union of independent Dawson-Watanabe processes with
initial conditions $\mu ,\nu$, respectively, and so the quantity on
the left side of \eqref{eq:weierstrass} is linear in $\mu$. It follows
that \eqref{eq:weierstrass} holds for all initial measures $\mu$ with
support contained in $D$.
\end{pf}

\section{Spatial Epidemic Models: Proof of Theorem
\ref{theorem:sis}}\label{sec:spatial}

\subsection{Strategy}\label{ssec:strategy} We shall exploit the fact
that the laws of the spatial \emph{SIS} and \emph{SIR} epidemics are
absolutely continuous with respect to the laws of their branching
envelopes.  Because the branching envelopes converge weakly, after
rescaling, to super-Brownian motions, by Watanabe's theorem (and in
the stronger sense of Theorem \ref{theorem:BRW}), to prove the weak
convergence of the rescaled spatial epidemics it will suffice to show
that the likelihood ratios converge, in a suitable sense, to the
likelihood ratios of the appropriate Dawson-Watanabe processes
relative to super-Brownian motion:

\begin{proposition}\label{proposition:strategy}
Let $X_{n}, X$ be random variables valued in a metric space
$\mathcal{X}$, all defined on a common probability space $(\Omega
,\mathcal{F},P)$, and let $L_{n},L$ be nonnegative, real-valued random 
variables  on $(\Omega ,\mathcal{F},P)$ such that 
\begin{equation}\label{eq:probabilityNormalization}
	E_{P}L_{n}=E_{P}L=1 \quad \forall \; n.
\end{equation}
Let $Q_{n},Q$ be the probability measures on $(\Omega ,\mathcal{F})$
with likelihood ratios $L_{n},L$ relative to $P$. If
\begin{equation}\label{eq:jointWeakConv}
	(X_{n},L_{n}) \Longrightarrow (X,L)
\end{equation}
under $P$ as $n \rightarrow \infty$, then the $Q_{n}-$distribution of
$X_{n}$ converges to the $Q-$distribution of $X$, that is, for every
bounded continuous function $f:\mathcal{X}\rightarrow \zz{R}$,
\begin{equation}\label{eq:qWeakConv}
	\lim_{n \rightarrow \infty}E_{Q_{n}}f (X_{n}) =E_{Q}f (X)
\end{equation}
\end{proposition}

\begin{proof}
Routine.
\end{proof}

Recall \cite{dawson:geo} that the law of the Dawson-Watanabe process
with killing is absolutely continuous with respect to that of the
standard Dawson-Watanabe process, with likelihood ratio given by
\eqref{eq:dwLR}.  The likelihood ratio involves stochastic integration
with respect to an orthogonal martingale measure; thus, the obvious
strategy for proving \eqref{eq:jointWeakConv} in our context is to
express the likelihood ratios of the spatial epidemic processes in
terms of stochastic integrals against the OMMs of the branching
envelopes.  Although it is possible to work directly with the
likelihood ratios of the epidemic processes to their branching
envelopes, this is somewhat messy, for two reasons: (A) the offspring
distributions $\mathcal{R}_{N}$ of the branching envelopes change with
the village size $N$; and (B) the random mechanism by which particles
of the branching envelope are culled in the standard coupling involves
dependent Bernoulli random variables. Therefore, we will first show,
by comparison arguments, that the spatial epidemic processes can be
modified so that difficulties (A) and (B) are circumvented, and in
such a way that asymptotic behavior is not affected. The likelihood
ratios of the modified processes relative to critical Poisson
branching random walks will then be computed in sec.~\ref{ssec:lr}.

\subsection{Extent, duration, size, and density of the branching
envelope}\label{ssec:duration} Since the spatial \emph{SIS} and
\emph{SIR} epidemics are stochastically dominated by their branching
envelopes, their durations, sizes,  etc., are limited by
those of their envelopes. For critical branching random walks, the
scaling limit theorems of Watanabe and Feller give precise information
about these quantities. Consider first the \emph{duration}: Since the
total mass in a \emph{BRW} is just a Galton-Watson process, if the
branching random walk is initiated by $N^{\alpha}$ particles, then by
a standard result in the theory of Galton-Watson processes (Th.~I.9.1
of \cite{athreya-ney}), the time $T_{N}$ to extinction scales like
$N^{\alpha}$, that is,
\begin{equation}\label{eq:GWduration}
	T_{N}/N^{\alpha} \Longrightarrow  F.
\end{equation}
The limit distribution $F$ is the distribution of the first passage
time to $0$ by a Feller diffusion process started at
$1$. Consequently, under the hypotheses of Theorem~\ref{theorem:sis},
the duration of the \emph{BRW} is $O_{P} (N^{\alpha})$. Furthermore,
if $Z^{N}_{n}$ is the mass in the $n$th
generation (that is, the total number of particles), then by Feller's
theorem (\cite{etheridge}, ch.~1), 
\begin{equation}\label{eq:fellerTheorem}
	Z^{N}_{[N^{\alpha}t]}/N^{\alpha}
	\Longrightarrow 
	Z_{t},
\end{equation}
where $Z_{t}$ is a Feller diffusion process started at
$1$. Consequently, the \emph{total mass} produced during the  entire
course of the branching envelope is of order $O_{P} (N^{2\alpha})$: in
particular, 
\begin{equation}\label{eq:totalMassScaling}
	\sum_{n}Z^{N}_{n}/N^{2\alpha}
	\Longrightarrow 
	\int_{0}^{\infty} Z_{t} \,dt.
\end{equation}
Since the Feller diffusion is absorbed at $0$ in finite time almost
surely, the integral is finite with probability $1$.

Next, consider the \emph{maximal density} and \emph{spatial extent} of
the branching random walk. Watanabe's theorem implies that if
initially all particles are located in an interval of size $N^{\alpha
/2}$ centered at $0$, as required by Theorem~\ref{theorem:sis}, then
the bulk of the mass must remain within $O_{P} (N^{\alpha/2})$ of the
origin, because the limiting Dawson-Watanabe process has bounded
support. A theorem of Kesten \cite{kesten:brw} implies that in fact
\emph{all} of the mass remains within $O_{P} (N^{\alpha/2})$ of the
origin: that is, under the hypotheses of Theorem~\ref{theorem:sis}, if
\begin{equation}\label{eq:maxBRW}
	Y^{*}_{N}:=
	\max \left\{|x|: \, \sum_{t<\infty }Y^{N}_{t} (x)>0\right\},
\end{equation}
then for any $\varepsilon >0$
there exists $\beta   <\infty$ such that 
\begin{equation}\label{eq:maxSpatialBRW}
	P\{ Y^{*}_{N}\geq  \beta  N^{\alpha /2}\} <\varepsilon .
\end{equation}
Together, \eqref{eq:GWduration} and \eqref{eq:maxBRW} imply that if
initially the branching random walk has $O (N^{\alpha})$ particles all
located at sites within distance $O (N^{\alpha /2})$ of the origin,
then the number of site/time pairs $(x,t)$ reached by the branching
random walk is $O_{P} (N^{3\alpha /2})$.  Now suppose in addition
that the initial configurations satisfy the more stringent
requirement $X^{N} (0,x)\Rightarrow X (0,x)$ of
Theorem~\ref{theorem:sis}; then Theorem~\ref{theorem:BRW} implies that
the renormalized density processes $X^{N} (t,x)\Rightarrow X (t,x)$,
where $X (t,x)$ is the renormalized density process of the standard
Dawson-Watanabe process. Since $X (t,x)$ is jointly continuous and has
compact support \cite{konno-shiga}, it follows that 
\begin{equation}\label{eq:maxDensity}
	\max_{t,x} Y^{N} (t,x)=O_{P} (N^{\alpha /2}).
\end{equation}

\subsection{Binomial/Poisson and Poisson/Poisson
comparisons}\label{ssec:BinomPoisson} In this section we show that, in
the asymptotic regimes considered in Theorem \ref{theorem:sis}, the
Binomial-$(N,1/N)$ random variables used in the the standard coupling
(sec.~\ref{ssec:heuristics}) can be replaced by Poisson-$1$ random
variables without changing the asymptotic behavior of the density
processes $X^{N} (t,x)$. Recall that in the standard coupling, each
particle, whether red or blue, produces a random number of offspring
with the Binomial-$(N,1/N)$ distribution. The total number of
particles produced during the lifetime of the branching envelope is,
under the hypotheses of Theorem \ref{theorem:sis}, at most $O_{P}
(N^{2\alpha})$, and $\alpha \leq 2/3$ in all scenarios considered, by
\eqref{eq:fellerTheorem}. Consequently, if all of the Binomial-$1/N$
random variables used in the construction were replaced by Poisson-$1$
random variables, the resulting processes (both the red process,
representing the spatial epidemic, and the red+blue process,
representing the branching envelope) would be indistinguishable from
the original processes, by the following lemma.

\begin{lemma}\label{lemma:Poisson-Binomial}
Assume that under the probability measure $\mu_{N}$, the random
variables $X_{1},X_{2},\dotsc$ are i.i.d.  Binomial-$(N,1/N)$, and
that under measure $\nu$ they are i.i.d. Poisson-$1$. Let $m=m_{N}$ be
a sequence of positive integers such that for some $\varepsilon >0$,
\begin{equation}\label{eq:mN}
	m_{N}=O (N^{2-\varepsilon})
\end{equation}
If $\mathcal{G}_{m}$ is the $\sigma -$algebra generated by
$X_{1},X_{2},\dotsc ,X_{m}$, then 
\begin{equation}\label{eq:BinomPoissonLR}
	\left(\frac{d\mu_{N}}{d\nu} \right)_{\mathcal{G}_{m}}
	\stackrel{\nu}{\longrightarrow}1. 
\end{equation}
\end{lemma}

\begin{proof}
This is a routine calculation. Fix a sequence
$x_{1},x_{2},\dotsc,x_{m}$ of nonnegative integers; the likelihood
ratio $d\mu_{N}/d\nu$ of this sequence is
\[
	\prod_{i=1}^{m}
	\frac{N!(N-1)^{-x_{i}}}{(N-x_{i})!} (1-N^{-1})^{N}e^{1}.
\]
For $m=O (N^{2-\varepsilon})$,
\[
	 \{ (1-N^{-1})^{N}e^{1} \}^{m}\sim \exp \{-m/2N \}.
\]
By Chebyshev's inequality and elementary calculus, for $\nu -$typical
sequences $x_{i}$,
\[
	\prod_{i=1}^{m}
	\frac{N!(N-1)^{-x_{i}}}{(N-x_{i})!}
	=\exp \{ (m+O_{P} (\sqrt{m}))/2N \}.
\]
\end{proof}

Recall that in the standard coupling (sec.~\ref{ssec:heuristics}), red
particles represent infections that occur in the spatial epidemic,
whereas blue particles represent attempted infections that are
suppressed because either two or more infected individuals try to
infect the same susceptible simultaneously, or (in the \emph{SIR}
case) because the target of the attempted infection has acquired
immunity by dint of an earlier infection. It is possible that more
than one attempted infection is suppressed at once, that is, more than
one blue particle with a red parent is created at a given site/time. In
sec.~\ref{ssec:collisions} below, we will show that such occurrences
are sufficiently rare that their effect on the epidemic process is
negligible in the large-$N$ limit.  To do so, we will bound the set of
offspring of such blue particles by the set of discrepancies between a
Poisson branching random walk with mean offspring number $1+\varepsilon$
and one with mean $1$, for some small $\varepsilon$. The next result
shows that if $\varepsilon$ is sufficiently small relative to the size
(total number of particles) of the Poisson branching random walk, then
the effect of changing the mean offspring number is negligible.

\begin{lemma}\label{lemma:Poisson-Poisson}
Let $\mu_{K}$ and $\nu_{K}$ be the distributions of Poisson branching
random walks with mean offspring numbers $1+\varepsilon_{K}$ and $1$,
respectively, and common initial configuration $Y^{K}_{0} (x)$ with
$K$ particles. If 
\begin{equation}\label{eq:meanDiscrepancy}
	\varepsilon_{K}=o (1/K)
\end{equation}
then under $\nu_{K}$, as $K \rightarrow \infty$.
\begin{equation}\label{eq:asymptoticallyIndistinuishable}
	\frac{d\mu_{K}}{d\nu_{K}} \Longrightarrow 1 .
\end{equation}
\end{lemma}

\begin{proof}
For a given sample evolution in which $Z_{K}$ particles are created, the
likelihood ratio of $\mu_{K}$ relative to $\nu_{K}$ is
\[
	\frac{d\mu_{K}}{d\nu_{K}}= 
	(1+\varepsilon_{K})^{Z_{K}}\exp \{-\varepsilon_{K}Z_{K} \}
	=1+O (Z_{K}\varepsilon_{K}^{2})
\]
Under $\nu_{K}$, the branching random walk will last on the order of
$K$ generations, during which on the order of $K$  particles will be
created in each generation, by Watanabe's theorem. Hence, $Z_{K}$ will
typically be of size $K^{2}$. In fact, if $\xi^{K}_{n}$ is the number of particle
creations in the $n$th generation, then $\xi^{K}_{n}$ is, under $\nu_{K}$, a
Galton-Watson process with offspring distribution  Poisson-$1$ and
initial condition $\xi^{K}_{0}=K$,  so  Feller's theorem implies
\[
	Z_{K}/K^{2}\Longrightarrow  \int_{0}^{\infty} \xi_{t}\,dt,
\]
where $\xi_{t}$ is a Feller diffusion process with initial state
$\xi_{0}=1$. The assertion now follows from the hypothesis
\eqref{eq:meanDiscrepancy}.
\end{proof}

The measures $\mu_{K},\nu_{K}$ can be coupled as follows: Start with
an initial configuration $Y^{K}_{0} (x)$, as in the lemma, and let
particles reproduce and move as in a branching random walk with
offspring distribution Poisson-$(1+\varepsilon_{K})$. Attach to each
particle $\zeta$ a Bernoulli-$\varepsilon_{K}/ (1+\varepsilon_{K})$
random variable $U_{\zeta}$. Assign colors \emph{green} or
\emph{orange} to particles according to the following rules: (A)
Offspring of green particles are always green. (B) An offspring
$\zeta$ of an orange particle is green if $U_{\zeta}=1$, otherwise is
orange. Then the process of orange particles evolves as a branching
random walk with offspring distribution Poisson-$1$, and the process
of all particles, green and orange, evolves as a \emph{BRW} with offspring
distribution Poisson-$(1+\varepsilon_{K})$. Denote by $Y^{K,G}_{t}(x)$
the number of \emph{green} particles at site $x$, time $t$, and by
$X^{K,G} (t,x)$ the renormalized density function obtained from
$Y^{K,G}_{t} (x)$ by the rule \eqref{eq:Yn}, with $k=K$.

\begin{corollary}\label{corollary:Poisson-Poisson}
Assume that the initial configurations $Y^{K}_{0} (\cdot)$ of the
branching random walks satisfy the hypotheses of Theorem
\ref{theorem:BRW}. If $\varepsilon_{K}=o (1/K)$, then as $k \rightarrow
\infty$, 
\begin{equation}\label{eq:maxDiscrepancy}
	\max_{t,x} X^{K,G} (t,x) \Longrightarrow 0.
\end{equation}
\end{corollary}

\begin{proof}
Denote by $X^{K} (t,x)$ the renormalized density process associated
with the branching random walk of \emph{orange} particles. By
Theorem~\ref{theorem:BRW}, the processes $X^{K}\Rightarrow X$ where
$X=X (t,x)$ is the density of a standard Dawson-Watanabe process. By
Lemma~\ref{lemma:Poisson-Poisson}, the density $X^{K}+X^{K,G}$ of the
orange+green particle branching random walk also converges to the
density of a standard Dawson-Watanabe process. Since there is only one
Dawson-Watanabe density process, it must be that $X^{K,G}$ converges
weakly to zero in sup norm.
\end{proof}

\subsection{Multiple collisions}\label{ssec:collisions} 

In the standard coupling (sec.~\ref{ssec:heuristics}) of a spatial
\emph{SIS} or \emph{SIR} epidemic with its branching envelope,
offspring of red particles at each site choose labels $j\in [N]$ at
random, which are then used to determine colors as follows: (A) If an
index $j$ is chosen by more than one particle, then all but one of
these are colored blue. (B) (\emph{SIR} model only) If index $j$ was
chosen at the same site in an earlier generation, then all particles
that choose $j$ are colored blue. We call events (A) or (B), where
offspring of red particles are colored blue, \emph{collisions}. At a
site/time where there are $\geq 2$ blue offspring of red particles
we say that a \emph{multiple collision} has occurred.  In this
section, we show that for either \emph{SIS} or \emph{SIR} epidemics,
up to the critical thresholds (see the statement of
Theorem~\ref{theorem:sis}), the effects of multiple collisions on the
evolution of the red particle-process are asymptotically
negligible. In particular, this will justify replacing the standard
coupling of sec.~\ref{ssec:heuristics} by the following modification,
in which at each time/site there is at most one blue offspring of a
red parent.

\medskip \noindent \textbf{Modified Standard Coupling:} Particles are
colored \emph{red} or \emph{blue}. Each particle produces a random
number of offspring, according to the Poisson-$1$ distribution, which
then randomly move either $+1,-1$, or $0$ steps from their birth
site. Once situated, these offspring are assigned colors according to
the following rules:
\begin{itemize}
\item [(A)] Offspring of blue particles are  blue; offspring of red
particles may be either red or blue. 
\item [(B)] At any site/time $(x,t)$ there is at most one blue
offspring of a red parent.
\item [(C)] Given that there are $y$ offspring of red parents at site
$x$, time $t$, the conditional probability
$\kappa_{N}(y)=\kappa_{N,t,x}(y)$ that one of them is blue is
\begin{align}\label{eq:kappaDefinitionSIS}
\kappa_{N}(y)&=y (y-1)/ (2N) \quad \text{for \emph{SIS}
epidemics}, \text{ and}\\
\label{eq:kappaDefinitionSIR}
\kappa_{N} (y)&=yR/N \quad \text{for \emph{SIR} epidemics}, \text{ where}\\
\label{eq:numberTagsUsed}
	R&=R^{N}_{t} (x) =\sum_{s<t}Y^{N}_{s} (x).
\end{align}
\end{itemize}
Here $Y^{N}_{t} (x)$ is the number of \emph{red} particles at site $x$
in generation $t$, and so $R=R^{N}_{t} (x)$ is the number of recovered
individuals at site $x$ at time $t$, equivalently, the number of
labels $j\in [N]$ that have been used in the standard coupling at $x$
by time $t$. Observe that in both the \emph{SIS} and the \emph{SIR}
cases, the value of $\kappa (y)$ is almost, but not exactly, equal to
the conditional probability that in the \emph{standard} coupling there
would be at least one blue offspring of a red particle. The small
discrepancies will make the expressions in the likelihood ratios
\eqref{eq:likelihoodratio} simpler.

The next result will justify replacing the standard coupling of
sec.~\ref{ssec:heuristics} by the modified standard coupling.
 
\begin{proposition}\label{proposition:collisions}
The standard couplings and the modified standard couplings can be
constructed simultaneously in such a way that the following is true,
for initial configurations satisfying the hypotheses of
Theorem~\ref{theorem:sis}.  If $Y^{N,\Delta}_{t} (x)$ is the
discrepancy between the numbers of red particles at $(x,t)$ in the
standard and modified standard couplings, then under the hypotheses of
Theorem \ref{theorem:sis}, as $N \rightarrow \infty$,
\begin{equation}\label{eq:collisionsAreNegligible}
	\max_{t,x} Y^{N,\Delta} (t,x) =o_{P} (N^{\alpha /2}).
\end{equation}
\end{proposition}

\begin{lemma}\label{lemma:4orMore}
Let $B_{N}$ ($B$ for ``bad'') be the number of sites/times $(x,t)$
where there are at least $4$ blue offspring of red particles in the
standard coupling. Then as
$N \rightarrow \infty$,
\begin{equation}\label{eq:4orMore}
	B_{N}=o_{P} (1).
\end{equation}
\end{lemma}

\begin{proof}
Theorem~\ref{theorem:BRW} implies that, under the hypotheses of
Theorem~\ref{theorem:sis}, the maximum number of particles (of any
color) at any site/time in the standard coupling is $O_{P}
(N^{\alpha/2})$. (Note: This also relies on the fact that the limiting
Dawson-Watanabe density process $X (t,x)$ is continuous and has
compact support, w.p.1.) Hence, we may restrict attention to sample
evolutions where $y=y_{t} (x)$, the total number of offspring of red
particles at $(x,t)$, satisfies $y\leq CN^{\alpha /2}$ for some fixed
constant $C<\infty$. Furthermore, by the considerations of
sec.~\ref{ssec:duration}, we may restrict attention to sample
evolutions of duration $\leq CN^{\alpha}$ and spatial extent
$CN^{\alpha /2}$. Since $\alpha \leq 2/3$, it follows that the number
of pairs $(x,t)$ visited by particles of the  branching envelope is no
more than $C^{2}N$.

In order that there be at least $4$ blue offspring of red parents at
site $x$, time $t$, at least 4 pairs (possibly overlapping) of
red-parent offspring must choose common labels $j\in [N]$.  The
conditional probability of this happening, given the value of
$y=y_{t}(x)$, is no more than $C'y^{8}/N^{4}$, for some constant $C'$
not depending on $y$ or $N$. But $y\leq CN^{\alpha /2}$, so this
conditional probability is bounded by $C''N^{4\alpha -4}\leq
C''N^{-4/3}$. Since there are only $C^{2}N$ sites to consider, it
follows that, on the event delineated in the preceding paragraph,  the
probability that $B_{N}\geq 1$ is $O (N^{-1/3})$.
\end{proof}

\begin{proof}
[Proof of Proposition~\ref{proposition:collisions}]
In this construction, each particle will be two-sided: the $S-$side
will represent the color of the  particle in the standard coupling,
and the $M-$side the color in the modified standard coupling. A
particle will be called a \emph{hybrid}  if the colors of its two
sides disagree. The strategy will be to show that colors can be
assigned in such a way that the process of hybrid particles is
dominated by the green particle process of
Corollary~\ref{corollary:Poisson-Poisson}; the result
\eqref{eq:collisionsAreNegligible} will then follow from
\eqref{eq:maxDiscrepancy}.

Consider first the \emph{SIS} case. Observe that in this case
$\kappa_{N} (y)$ is the conditional expectation, in the standard
coupling, of the number of pairs that share labels, given that there
are $y$ offspring of red parents at a site. Thus, $\kappa_{N} (y)$
exceeds (by a small amount) the conditional probability in the
standard coupling that at least one of the $y$ red-parent offspring
would be blue. Denote by $\Delta_{N} (y)$ the excess; note that
\[
	\Delta_{N} (y)=O (y^{4}/N^{2}).
\]
 
The rules by which the process evolves are as follows: First, all
particles reproduce, each creating a random number of offspring with
Poisson-$1$ distribution.  Each offspring then moves $+1,-1$, or $0$
steps from its birth site, and chooses a random label $j\in
[N]$. Particles with ``genotype'' $BB$ (that is, offspring of
particles with coloring $BB$; the first letter denotes the $S-$color,
the second the $M-$color) will always be colored $BB$, and their
labels $j$ will play no role in determining the colors of the other
offspring. However, the labels of all other offspring matter. Say that
there is an $S-$\emph{duplication} at label $j$ if at least two
particles both with a red $S-$gene choose label $j$; similarly, say
that an $M-$\emph{duplication} occurs at $j$ if $j$ is chosen by at
least two particles with $M-$gene $R$. (Note: If both $BR$ and $RB$
genotype particles choose label $j$, only the $RB$ particles are
counted  in the possible $S-$duplication, and only $BR$ particles in the
$M-$duplication.) Particles at $(x,t)$ are now
assigned color ``phenotypes'' by the following rules:
\begin{itemize}
\item [(D)] (Default) If there is a duplication involving at least two
particles of genotype $RR$, do the following: Among all such
duplications, choose one (say $i$) at random; choose one of the
genotype-$RR$ particles with label $i$, assign it phenotype $BB$, and
give all of the other particles with label $i$ the same genotypes as
their parents. Give all other particles at the site the same
$M-$colors as their parents, and assign the remaining $S-$colors by
rule (S):
\item [(S)] If there are labels $j\not =i$ with $S-$duplications, then for
each such label 
\begin{itemize}
\item [(a)] If there is at least one particle with genotype $RB$ in
the duplication, then choose one of all such $RB$ particles at random, 
give it phenotype $BB$, and give all other particles involved in the
duplication $S-$color $R$.
\item [(b)] Otherwise, if all particles involved in the duplication
have genotype $RR$, choose one at random and assign it phenotype $BR$,
and give all of the rest phenotype  $RR$. 
\item [(c)] Give all particles not involved in $S-$duplications
the same $S-$colors as their parents.
\end{itemize}
\item [(M)] If there are no duplications of type (D) but at least one
$M-$duplication, then in any such duplication, at least one particle
of genotype $BR$ must be involved. Choose one at random and
assign it phenotype $BB$,   give all of the remaining particles at
the site $M-$color $R$, and assign $S-$colors by rules (S)-(a),(c).
\item [(A)] (Adjustment Step) If there are no $M-$duplications, then
toss a $\Delta_{N} (y)-$coin: If it comes up Heads, choose one of the
particles with genotype $?R$ at random, give it $M-$color $B$, and
give all of the rest $M-$color $R$. If it comes up Tails, give every
particle the same $M-$color as its parent.
\end{itemize}

These rules guarantee that the $S-$colors of the particles behave as
in the standard coupling, and that the $M-$colors behave as in the
modified standard coupling. Therefore, the discrepancy
$Y^{N,\Delta}_{t} (x)$ is bounded by the number of hybrid particles at
$(x,t)$. Hybrids can be offspring of $RR$, $BR$, or $RB$ particles,
but not $BB$ particles; however, a hybrid can only be produced by a
particle of type $RR$ if (i) there is a multiple collision, i.e., if
there are at least two pairs of non-$BB$ particles that choose the
same labels; or (ii) the coin toss in the adjustment step (A) is a
Head. Both of these are events of (conditional) probability
$O(y^{4}/N^{2})$.  Moreover, by Lemma~\ref{lemma:4orMore}, except with
vanishingly small probability, there is no site/time $(x,t)$ with more
than 3 hybrid offspring of $RR$ parents. Consequently, on the event 
that the maximal number of particles at any site/time is no more than 
$CN^{\alpha /2}$ (see \eqref{eq:maxDensity} above), the process of
$RR \rightarrow$hybrid creations is dominated as follows: Let each particle,
in every generation, produce an additional Poisson-$2CN^{3\alpha/2-2}$ offspring;
immediately replace each such particle by 3 \emph{green} particles, and let
green particles only beget other green particles in subsequent
generations. Since $\alpha \leq 2/3$, the rate at which green
particles are produced by non-green particles is $O (N^{-1})$, so
Corollary~\ref{corollary:Poisson-Poisson} implies that the green
particle process is asymptotically negligible. This proves
\eqref{eq:collisionsAreNegligible}  in the \emph{SIS} case. The
\emph{SIR} case is proved by a very similar construction. 
\end{proof}

The upshot of Proposition~\ref{proposition:collisions} is that in
proving Theorem~\ref{theorem:sis}, the epidemic process (the red
process in the extended standard coupling) can be
replaced by  the red process in the modified standard coupling.

\subsection{Orthogonal martingale measures and convergence of
stochastic integrals}\label{ssec:omm}

Let $Y^{k}_{t} (x)$ be the number of particles at site $x\in \zz{Z}$
at time $[t]\in \zz{Z}_{+}$ in a one-dimensional branching random walk
with Poisson-$1$ offspring law $\mathcal{R}_{\infty}$. Denote by
$X^{k}_{t}= \mathcal{F}_{k}Y^{k}_{kt}$ the corresponding rescaled
measure-valued process.  For each $k$, the measure-valued process
$X^{k}_{t}$ satisfies a martingale problem analogous to that satisfied
by the super-Brownian motion: If $\phi \in C^{\infty}_{b}$, the
(cadlag) process
\begin{equation}\label{eq:discreteMGProblem}
	M^{k}_{t} (\phi):=
	\xbrace{X^{k}_{t},\phi}-\xbrace{X^{k}_{0},\phi} -
	\int_{0}^{[kt]/k}\xbrace{X^{k}_{s},{A}_{k}\phi } \, ds
\end{equation}
is a martingale, where $A_{k}$ is the difference operator
\begin{equation}\label{eq:Ak}
	A_{k}\phi (x)=\{\phi (x+1/\sqrt{k})+\phi (x-1/\sqrt{k})-2\phi
	(x))\}/{3k^{-1}}. 
\end{equation}
(Since  $X_{t}$ is constant on successive time intervals of
duration $k^{-1}$, the integral in \eqref{eq:discreteMGProblem} is really a
sum.)  The operator $\phi \mapsto M^{k}_{.}
(\phi)$ extends to an \emph{orthogonal martingale measure} $M^{k}
(ds,dx)$ (see \cite{walsh} for the definition and basic stochastic
integration theory). The
measure $M^{k}$ is purely discrete, putting mass only at points
$(s,x)\in k^{-1}\zz{Z}_{+}\times k^{-1/2}\zz{Z}$: at such points
$s=n/k,x=m/\sqrt{k}$, 
\begin{align}\label{eq:ommExplicit}
	{k}M^{k} (\{s \}, \{x \})
	 &=Y^{k}_{n} (m)-\lambda^{k}_{n} (m)
	 \quad 
	 \text{where}\\
\label{eq:lambdaDef}
      \lambda^{k}_{n} (m)&=\lambda^{k} (s,x):
      =E (Y^{k}_{n} (m)\, | \, \mathcal{H}_{n-1})
      =\frac{1}{3}\sum_{i=0,-1,1}Y^{k}_{n-1} (x+i), \quad \text{and}\\
\label{eq:filtDef}
	 \mathcal{H}_{n}&=\sigma (\{ Y^{k}_{j} (m)\, : \, m\in \zz{Z},j\leq n\})
\end{align}
is the $\sigma -$algebra generated by the history of the evolution to
time $n$. Note that mass is scaled by the factor $k$, as required
by the Feller-Watanabe normalization. Note also that, conditional on
$\mathcal{H}_{n-1}$, the random variables $\{ Y^{k}_{n} (m)\}_{m\in
\zz{Z}}$ are mutually independent: this implies that the martingale
measure $M^{k}$ is an \emph{orthogonal} martingale measure.

\begin{proposition}\label{proposition:OMMConvergence}
Assume that the initial particle densities $X^{k} (0,\cdot)$ satisfy
the hypotheses of Theorem \ref{theorem:BRW}, that is, they have common
compact support and they converge in $C_{b} (\zz{R})$ to a continuous
function $X (0,\cdot)$. Then the random vectors $(X^{k},M^{k})$
consisting of the density functions $X^{k} (t,x)$ and the orthogonal
martingale measures $M^{k}$ converge weakly as $k \rightarrow \infty$
to $(X,M)$, where $X$ is the Dawson-Watanabe density process
(super-Brownian motion) with initial condition $X (0,x)$, and $M$ is
its associated orthogonal martingale measure.
\end{proposition}

\begin{proof}
Consider first the marginal distributions of the orthogonal martingale
measures $M^{k}$, viewed as random elements of the Skorohod space
$\zz{D} ([0,\infty ),\mathcal{S}')$, where $\mathcal{S}'$ is the space
of tempered distributions on $\zz{R}$. In order to prove that $M^{k}
\Rightarrow M$ it suffices, by Mitoma's theorem (cf. \cite{walsh},
Th.~6.15), to prove that (i) for any $\phi\in \mathcal{S}$ (where
$\mathcal{S}$ is the Schwartz space of test functions), the processes
$M^{k}_{t} (\phi )$ are tight, and (ii) finite-dimensional
distributions $M^{k}_{t_{i}} (\phi_{i})$ converge for all $\phi_{i}\in
\mathcal{S}$. Both of these follow routinely from the representation
\eqref{eq:discreteMGProblem} and Watanabe's theorem: In particular,
Watanabe's theorem implies that for any finite subset $\{\phi_{i} \}_{i\in I}$
of $\mathcal{S}$, 
\[
	(\xbrace{X^{k}_{t},\phi_{i}}-\xbrace{X^{k}_{0},\phi_{i}})_{i\in
	I} \Longrightarrow 
	(\xbrace{X_{t},\phi_{i}}-\xbrace{X_{0},\phi_{i}})_{i\in I}
\]
and
\[
	\left(\int_{0}^{t}\xbrace{X^{k}_{s},A\phi_{i}} \,ds\right)_{i\in I}
	\Longrightarrow 
	\left(\int_{0}^{t}\xbrace{X_{s},\phi_{i}}\,ds \right)_{i\in I}
\]
where
\[
	A=\lim_{k \rightarrow \infty}A_{k}=\Delta /\sqrt{6}.
\]
Consequently, to deduce (i)--(ii) above it suffices to show that for
each $\phi =\phi_{i}$,
\[
	\int_{0}^{t} \xbrace{X^{k}_{s},|A_{k}\phi-A\phi | }\, ds
	+
	\int_{[kt]/k}^{t}\xbrace{X^{k}_{s},|A\phi |} \, ds
	\stackrel{P}{\longrightarrow}0. 
\]
Since $\xnorm{A_{k}\phi-A\phi}_{\infty}\rightarrow 0$ for
Schwartz-class functions $\phi$, this also follows from Watanabe's
theorem.

It remains to show that the convergence $M^{k}\Rightarrow M$ holds
jointly with $X^{k} (t,x)\Rightarrow X (t,x)$.  Since
$X^{k}\Rightarrow X$ (Theorem~\ref{theorem:BRW}) and
$M^{k}\Rightarrow M$ marginally, the joint distributions are
tight. Hence, to prove that $( X^{k},M^{k})$ converge \emph{jointly}, it
suffices to show that the only possible weak limit is $(X,M)$, where
$X=X (t,x)$ is the density process associated with the Dawson-Watanabe
process. But because a continuous function $w (x)$  is determined by
its integrals against Schwartz-class functions $\phi$, it suffices to
show that finite dimensional distributions of the vector-valued
processes $( M^{k}_{t}\phi_{i},\xbrace{X^{k}_{t},\phi_{i} })_{i \in I}$ converge
to the  corresponding  joint distributions of
$(M_{t}\phi_{i}, \xbrace{X_{t},\phi_{i}})_{i\in I}$. This follows by a
repetition of the  argument in the preceding paragraph.
\end{proof}

\begin{corollary}\label{corollary:SIConvergence}
Let $\theta (t,x,w)$ be a bounded, jointly continuous function of
$t\geq 0$, $x\in \zz{R}$, and $w\in \zz{D}([0,0),C_{c} (\zz{R}))$ such
that for any $t$ the function $\theta (t,x,w)$ depends only on
$w[0,t]$, that is, $\theta (t,x,w)=\theta (t,x,w[0,t])$.  Assume that
the hypotheses of Proposition \ref{proposition:OMMConvergence}
hold. Then
\begin{equation}\label{eq:SIA}
	\iint \theta (s,x, X^{k}) \,M^{k}
	(ds,dx)\Longrightarrow 
	\iint \theta (s,x,X)\, M (ds,dx) 
\end{equation}
Moreover,  \eqref{eq:SIA} holds jointly with the convergence
$X^{k}\Rightarrow X$ in $\zz{D} ([0,\infty ),C_{b} (\zz{R}))$.
\end{corollary}

\begin{proof}
This can be deduced from Prop.~7.6 of \cite{walsh}, but verification
of the hypothesis (7.5) is more work than a direct proof. The elementary
Prop.~7.5 of \cite{walsh} implies that weak convergence \eqref{eq:SIA}
(and joint convergence with $X^{k}\Rightarrow X$) holds for
\emph{simple} integrands
\[
	\theta (t,x,w)=\sum_{i=1}^{n}
	       a_{i} (w)\mathbf{1}_{(s_{i},t_{i}]} (t)\phi_{i} (x)
\]
such that (a) each $\phi_{i}\in \mathcal{S}$, (b) each $a_{i} (w)$ is
bounded, continuous in $w$, and $\zz{F}_{s_{i}}-$measurable (here
$(\zz{F}_{s})$ is the natural filtration on $\zz{D}([0,t],C_{c}
(\zz{R}))$), and (c) none of the jump times $s_{i},t_{i}$ coincides
with jumps of one of the martingale measures $M^{k}$. Clearly, any
function $\theta$ satisfying the hypotheses of
Corollary~\ref{corollary:SIConvergence} can be \emph{uniformly}
approximated by such simple functions. Consequently, to prove
\eqref{eq:SIA} it suffices to show that for any $\varepsilon >0$ there
exists $\delta =\delta (\varepsilon)>0$ such that for any simple
function $\theta$ satisfying (a)-(c) and $\xnorm{\theta}_{\infty}<\delta$,
\begin{equation}\label{eq:SISmall}
	P\left\{\bigg| \int \theta (t,x)M^{k} (dt,dx)\bigg|>
	\varepsilon\right\}
	 <\varepsilon \quad \forall \,k.
\end{equation}

To establish \eqref{eq:SISmall} we use the special structure of the
orthogonal martingale measure $M^{k}$. For each $k$, this is a purely
discrete random measure with atoms \eqref{eq:ommExplicit}. By
hypothesis, the conditional distribution of $Y^{k}_{n} (m)$ given the
past is Poisson with mean (and therefore also variance)
$\lambda^{k}_{n} (m)$ (see \eqref{eq:lambdaDef}), and the random
variables $\{ Y^{k}_{n} (m)\}_{m\in \zz{Z}}$ are, for each fixed $n$,
conditionally independent given the past. Hence, the predictable
quadratic variation of the local martingale
\[
	(\theta \cdot M^{k})_{t}:=\int_{0}^{t}\int_{x}\theta
	(s,x)M^{k} (ds, dx)
\]
is 
\[
	[\theta \cdot M^{k}]_{t}=
	\sum \sum \mathbf{1}_{[0,t)} (s)\theta (s,x)^{2} \lambda^{k} (s,x)
	     /k^{2},
\]
where the sum is over the jump points $(s,x)\in
k^{-1}\zz{Z}_{+}\times k^{-1/2}\zz{Z}$ of the martingale measure
$M^{k}$. Thus, if $\xnorm{\theta}_{\infty}<\delta$, then the quadratic
variation of $(\theta \cdot M^{k})_{t}$ is bounded by
\[
	\delta^{2}\sum \sum \lambda^{k} (s,x)/k^{2}
	=\delta^{2}\sum \sum Y^{k} (s,x)/k^{2}
\]
But $\sum \sum  Y^{k} (s,x)/k^{2}$ is just the total rescaled mass in
the branching random walk, which by Watanabe's theorem (or  Feller's
theorem) converges in law to the total mass in the  standard
Dawson-Watanabe process. The inequality \eqref{eq:SISmall} now follows
routinely.
\end{proof}

Unfortunately, the functions for which we would like to apply this
result --- namely, those defined by equations
\eqref{eq:sirKillingRate} and \eqref{eq:sisKillingRate} -- are
unbounded. Worse, the function $\theta (t,x)=X (t,x)$ in equation
\eqref{eq:sisKillingRate}  isn't even continuous (as a function of
$X(t,x)\in \zz{D} ([0,\infty )\times C_{c} (\zz{R}))$, because such
functions may have jumps). The following corollary takes care of the
first problem.

\begin{corollary}\label{corollary:unboundedSIConvergence}
Assume that the hypotheses of
Proposition~\ref{proposition:OMMConvergence} are satisfied.
Let $\theta (t,x,w)=\theta (t,x,w[0,t])$ be a  jointly continuous
function of $t\geq 0$, $x\in \zz{R}$, and $(w\in \zz{D}
([0,\infty ),C_{c} (\zz{R}))$ such that  for every scalar $C<\infty$ and every
compact subset $F$ of $[0,\infty )\times \zz{R}$,
\begin{equation}\label{eq:uniformBoundedness}
	\sup \{|\theta (t,x,w)| \, : \,  \sup_{t',x'}|w (t',x')|\leq C
	     \; \text{and \rm support} (w)\subset F\; \}
	     <\infty.
\end{equation}
Then \eqref{eq:SIA} holds  jointly with the convergence
$X^{k}\Rightarrow X$ in $\zz{D} ([0,\infty ),C_{b} (\zz{R}))$.
\end{corollary}

Observe that the hypothesis  is satisfied by the function $\theta$
defined by \eqref{eq:sirKillingRate}.

\begin{proof}
Because the Dawson-Watanabe density process $X (s,x)$ is almost surely
continuous with compact support, continuity of $\theta$ ensures that
\[
	\iint  \theta (s,x,X)^{2}X (s,x) \,
	ds \,dx <\infty  ,
\]
and this in turn guarantees that
\begin{equation}\label{eq:SITruncation}
	\lim_{C \rightarrow \infty}
	\iint (\theta (s,x,X)\wedge C) \,M
	(ds,dx)
	= \iint  \theta (s,x,X) \,M (ds,dx)
\end{equation}
and that the limit is a.s. finite.  By the hypothesis \eqref{eq:uniformBoundedness},
the function $\theta$ is uniformly bounded on any set of sample
evolutions $\omega =(x_{s})$ such that support$(\omega)$ and $\sup
|\omega |$ are uniformly bounded. By the results of
sec.~\ref{ssec:duration} above, the supports and suprema of the random
functions $X^{k} (t,x)$ are tight, that is, for any $\varepsilon >0$
there exist a compact set $F=F_{\varepsilon}$ and a constant
$C=C_{\varepsilon}<\infty$ such that for all $k$,
\[
	 P\{\text{support} (X^{k})\subset F \; \text{and}\; \max X^{k}
	 (t,x)\leq C\} \geq 1-\varepsilon .
\]
 Consequently, since the support of the martingale measure $M^{k}$ is
 contained in that of $X^{k}$, it follows that for any
$\varepsilon>0$ there exists $C<\infty$ such that  for all $k$,
\[
	 \iint  \theta (s,x,X) \,M^{k} (ds,dx)=
	 \iint  (\theta (s,x,X)\wedge C) \,M^{k} (ds,dx)
\]
except on a set of probability $\leq \varepsilon$.  Weak convergence
\eqref{eq:SIA} (and joint weak convergence with $X^{k} \Rightarrow X$)
now follows routinely from Corollary~\ref{corollary:SIConvergence} and
\eqref{eq:SITruncation}. 
\end{proof}

The function $\theta (t,x)=X (t,x)/2$ that occurs in the \emph{SIS}
case of Theorem~\ref{theorem:sis} is not continuous for functions
$X(t,x)$ in the  space $\zz{D} ([0, \infty ), C_{c} (\zz{R}))$,
because such functions may have jumps. Thus,
Corollaries~\ref{corollary:SIConvergence}--\ref{corollary:unboundedSIConvergence}
do not apply directly. The following corollary addresses this case
specifically.

\begin{corollary}\label{corollary:SIConvergenceSISCase}
Assume that the hypotheses of
Proposition~\ref{proposition:OMMConvergence} are satisfied, and define
$\theta^{k} (s,x)=\lambda^{k} (s,x)/\sqrt{k}$, where
$\lambda^{k}(s,x)$ is as in \eqref{eq:lambdaDef}  above. Then 
\begin{equation}\label{eq:SIConvergenceSISCase}
	\iint \theta^{k} (s,x) \,M^{k} (ds,dx)
	\Longrightarrow 
	\iint X (s,x)\,M (ds,dx),
\end{equation}
and this convergence holds jointly with $X^{k}\Rightarrow X$.
\end{corollary}

\begin{proof}
Although the random functions $X^{k} (t,x)$ and $\lambda^{k} (t,x)$
have jumps, the jumps are, with high probability, small, because
$X^{k} (t,x)\Rightarrow X (t,x)$, by
Theorem~\ref{theorem:BRW}. Consequently, it is still possible to
approximate $\theta^{k}$ by bounded, continuous functions $\varphi$ in
such a way that the stochastic integrals of $\theta^{k}$ relative to
$M^{k}$ are well approximated, with high probability, by those of
$\varphi$. In particular, define, for any $\varepsilon >0$ and
$C<\infty$,
\[
	\varphi (s,x)=\varphi_{C,\varepsilon} (s,x,w)
	=\frac{1}{2\varepsilon^{2}}
	\iint_{\substack{s'\in [s-2\varepsilon ,s-\varepsilon]\\
			       x'\in [x-\varepsilon ,x+\varepsilon]}}
			       w (s',x')\wedge C \,ds'\,dx';
\]
then for any $\delta >0$ there exist
$C,\varepsilon$ such that for all $k$,
\begin{equation}\label{eq:integrandApprox}
	P\{\max_{s,x}|\theta^{k} (s,x)- \varphi (s,x,X^{k}) |>\delta\}<\delta .
\end{equation}
Clearly, $\varphi (s,x,w)$ is jointly continuous in its arguments and
uniformly bounded by $C$, so it meets the requirements of
Corollary~\ref{corollary:SIConvergence}. It follows that for any
$C,\varepsilon$, 
\[
	\iint \varphi (s,x,X^{k})\,M^{k} (ds,dx)
	\Longrightarrow 
	\iint \varphi (s,x,X))\, M (ds,dx),
\]
and this holds jointly with $X^{k}\Rightarrow X$. Thus, to prove the corollary,
it suffices to show that if $C$ and $\varepsilon$ are suitably chosen
then the differences between the stochastic integrals of $\theta^{k}$
and $\varphi$ against the martingale measures $M^{k}$ are small with
high probability, uniformly in $k$.

By virtually  the same calculation as in the proof of
Corollary~\ref{corollary:unboundedSIConvergence}, the local martingale 
\[
	\iint_{s\leq t} (\theta^{k} (s,x)-\varphi (s,x,X^{k}))\,M^{k} (ds,dx)
	:=\iint_{s\leq t}\Delta^{k} (s,x)\, M^{k} (ds,dx)= (\Delta^{k}\cdot M^{k})_{t}
\]
has predictable quadratic variation
\begin{align*}
	[\Delta^{k}\cdot M^{k}]_{t}&=
			 \sum \sum  \mathbf{1}_{[0,t)} (s) \Delta^{k}
			 (s,x)^{2} \lambda^{k} (s,x)/2		\\
			 &\leq 
			 \max_{s,x}\Delta^{k} (s,x)^{2}\sum \sum 
			 Y^{k} (s,x)/k^{2}.
\end{align*}
If $\varepsilon >0$ is sufficiently small and $C<\infty$ sufficiently
large, then by inequality \eqref{eq:integrandApprox},
$\xnorm{\Delta^{k}}_{\infty}<\delta$ except with small probability,
uniformly for all $k$.  Recall that the sum $\sum \sum
Y^{k}(s,x)/k^{2}$ is the total rescaled mass in the branching process,
and so by Feller's theorem converges in law to the total mass in the
standard Dawson-Watanabe process. Hence, with high probability, the
quadratic variation $[\Delta^{k}\cdot M^{k}]_{\infty}$ will be small,
provided $C$ and $\varepsilon $ are chosen so that $\delta$ is small.
Therefore, by standard martingale arguments, the maximum modulus of
the stochastic integral $(\Delta^{k}\cdot M^{k})_{t}$ will be small,
with high probability, uniformly in $k$.
\end{proof}

\subsection{Likelihood ratios: generalities}\label{ssec:lr} The
strategy of the proof of Theorem~\ref{theorem:sis} is to show that the
likelihood ratios of the (modified) spatial epidemic processes
relative to their branching envelopes converge weakly to the
likelihood ratios of the appropriate Dawson-Watanabe processes with
killing relative to the Dawson-Watanabe process with no killing. The
likelihood ratio of the Dawson-Watanabe process with killing is given
by \eqref{eq:dwLR}. This expression involves a stochastic integral
relative to the orthogonal martingale measure of a standard
Dawson-Watanabe process, and so to prove weak convergence we will
express the likelihood ratios of the epidemic processes in terms of
stochastic integrals.

Consider a sequence $Y^{N}_{t} (x)$ of counting processes, and
probability measures $P=P^{N},Q=Q^{N}$ such that under $P$ the process
$Y^{N}$ is a branching random walk with offspring law
$\mathcal{R}_{\infty}$, and under $Q$ it is a modified epidemic
process (that is, the red particle process in the modified standard
coupling of sec.~\ref{ssec:collisions}). Assume that the initial
conditions $Y^{N}_{0}$ are common under $P$ and $Q$, and satisfy the
hypotheses of Theorem~\ref{theorem:sis}. Let
$\mathcal{H}_{t}$ be the $\sigma -$algebra generated
by the history of the evolution to time $t$, and set
\begin{align}\label{eq:lambda}
      \lambda^{N}_{t} (x):&= E_{P} (Y^{N}_{t} (x)\, | \, \mathcal{H}_{t-1})\\
 \notag   		             &=(Y^{N}_{t-1}(x-1)+Y^{N}_{t-1} (x)
		 	     		   +Y^{N}_{t-1} (x+1) )/3.
\end{align}

Under $P$, the random variables $Y^{N}_{t} (x)$ are, for each $t$,
conditionally independent Poisson r.v.s, given $\mathcal{H}_{t-1}$,
with conditional means $\lambda^{N}_{t} (x)$. In the modified standard
couplings, the color choices at the various sites $x$ are, conditional
on $\mathcal{H}_{t-1}$ and on the numbers of red-parent offspring at
the various sites, mutually independent, with at most one blue
offspring of a red parent at any site/time. Hence, the event
$Y^{N}_{t} (x)=y$ could occur in one of only two ways: (1) there are
$y$ offspring of red parents, none of which takes the color blue; or
(2) there are $y+1$ offspring of red parents, one of which is blue. It
follows that the relative likelihood $L_{N}:= dQ^{N}/dP^{N}$ of a
sample evolution $(y_{t} (x))_{t,x}$ is given by
\begin{equation}\label{eq:likelihoodratio}
	L_{N}=\prod_{t\geq 1}\prod_{x\in \zz{Z}}L_{N} (t,x)
			    :=\prod_{t\geq 1}\prod_{x\in \zz{Z}}
			    \frac{p (y|\lambda ) (1-\kappa_{N} (y))
			     +p ((y+1)|\lambda)\kappa_{N} (y+1)}{p (y|\lambda )}
\end{equation}
where in each factor, $y=y_{t} (x)$ and  $\lambda =\lambda^{N}_{t} (x)$,
and $\kappa_{N} (y)=\kappa_{N,t,x} (y)$ is the conditional probability
in the modified coupling, given $y$ offspring of red parents at
$(t,x)$, that one of these is colored blue. Here $p (\cdot |\lambda )$
is the Poisson density with parameter $\lambda$. The conditional
probability $\kappa_{N}(y)$ is given by \eqref{eq:kappaDefinitionSIS}
for \emph{SIS} epidemics, and by \eqref{eq:kappaDefinitionSIR} for
\emph{SIR} epidemics.

Although the product \eqref{eq:likelihoodratio} extends over
infinitely many sites and time x, t, all but finitely many of the
factors are 1: In particular, if the nearest neighbors (in the
previous generation) $(x',t-1)$ of site $(x,t)$ are devoid of
particles, then the factor $L_{N} (t,x)$ indexed by $(x,t)$ must be
$1$. Recall (by \eqref{eq:GWduration} and \eqref{eq:maxSpatialBRW} of
sec.~\ref{ssec:duration}) that under the hypotheses of
Theorem~\ref{theorem:sis}, the number of site/time
pairs $(x,t)$ visited by the branching envelope is of order
$O_{P} (N^{3\alpha /2})$. Thus, the number of nontrivial factors in
\eqref{eq:likelihoodratio}  is typically of the same order.

\subsection{Likelihood ratios: \emph{SIR} epidemics}\label{ssec:lrSIR}
Analysis of the likelihood ratios is somewhat easier in the \emph{SIR}
case, in that the error terms are more easily disposed. This is
because in the \emph{SIR} case Theorem~\ref{theorem:sis} requires that
$\alpha \leq 2/5$, and so the number of sites that contribute to  the
product \eqref{eq:likelihoodratio} is $O_{P} (N^{3/5})$. In particular,
error terms of order $o_{P} (N^{-3/5})$ can be ignored.

Recall that in the modified standard coupling for the \emph{SIR}
epidemic, the conditional probability
$\kappa_{N}(y)=\kappa_{N,t,x}(y)$ that there is a blue offspring of a
red parent at $(t,x)$ is $yR/N$, where $R=R^{N}_{t} (x)$ is the number
of labels used previously at site $x$ (see \eqref{eq:numberTagsUsed}).
Hence,  the contribution
to the likelihood ratio from the site $(t,x)$, on the sample evolution
$Y^{N}_{t} (x)=y_{t} (x)$, is
\begin{align}\label{eq:lrSIR}
L_{N} (t,x)&=1	-yR/N + (\lambda / (y+1)) (y+1)R/N\\
\notag 
&=1- (y-\lambda)R/N.
\end{align}
Thus, in the \emph{SIR}
case, the likelihood ratio \eqref{eq:likelihoodratio} can be written as
\begin{align}\label{eq:LN}
	L_{N}&= \prod_{t}\prod_{x}
	       (1- (( Y^{N}-\lambda^{N}_{t}(x))R^{N}_{t} (x)/N)\\
\notag  &=\prod_{t}\prod_{x}
	     \exp \{ \log (1-\Delta^{N}_{t}(x)\varrho^{N}_{t} (x))\}\\
\notag &=(1+\varepsilon_{N})
		\exp \left\{ -\sum_{t}\sum_{x}\Delta^{N}_{t}(x) \varrho^{N}_{t} (x)
		 -\frac{1}{2} \sum_{t}\sum_{x}\Delta^{N}_{t}(x)^{2}\varrho^{N}_{t} (x)^{2}
		\right\},
\end{align}
where 
\begin{align}\label{eq:DeltaN}
	\Delta^{N}_{t} (x)&=( Y^{N}_{t} (x) -\lambda^{N}_{t} (x))/N^{\alpha },\\
\label{eq:rhoN}
	\varrho^{N}_{t} (x)&=R^{N}_{t} (x)/N^{1-\alpha},\quad \text{and}\\
\label{eq:epsilonN}
	\varepsilon_{N}&=o_{P} (1).
\end{align}
(Note: The error in the two-term Taylor series approximation to the
logarithm is of magnitude $O_{P} (Y^{3}R^{3}/N^{3})=O_{P}
(N^{9\alpha/2-3})$, which is asymptotically negligible for $\alpha
\leq 2/5$; hence \eqref{eq:epsilonN}.)

The two sums that occur in the last exponential in \eqref{eq:LN} are a
stochastic integal and its corresponding quadratic variation,
respectively. To see this, observe that the quantities $\Delta^{N}_{t}
(x)$ coincide with the masses \eqref{eq:ommExplicit} in the orthogonal
martingale measures $M^{N}$ associated with the branching random walks
$Y^{N}$. (In \eqref{eq:ommExplicit}, $k=N^{\alpha}$; it makes more
sense here to index by $N$ rather than $k$.) Consequently,
\begin{equation}\label{eq:discreteSI}
	\sum_{t}\sum_{x}\Delta^{N}_{t} (x) \varrho ^{N}_{t} (x)
	=\iint \theta^{N}(t,x) \, M^{N} (dt,dx).
\end{equation}
where 
\begin{equation}\label{eq:thetaN}
	\theta^{N} (t,x)=\varrho^{N}_{tN^{\alpha }} (xN^{\alpha /2})
\notag 	       =N^{5\alpha /2-1}\int_{s=0}^{t}X^{N} (s,x)\,ds 
\end{equation}
For $\alpha <2/5$, this converges to zero as $N \rightarrow \infty$;
for $\alpha =2/5$ it coincides with $\theta (t,x)=\theta (t,x,X^{N})$
as given by \eqref{eq:sirKillingRate}. Therefore,
Corollary~\ref{corollary:unboundedSIConvergence}  implies that
\begin{equation}\label{eq:Cconvergence}
	\iint \theta^{N} (t,x)\, M^{N} (dt,dx)
	\Longrightarrow 
	\iint \theta (t,x) \, M (dt, dx),
\end{equation}
where $M$ is the orthogonal martingale measure attached to the
standard Dawson-Watanabe process and $\theta (t,x)$ is as in parts
(c)-(d) of
Theorem~\ref{theorem:sis}. (Corollary~\ref{corollary:unboundedSIConvergence}
also implies that the convergence holds jointly with $X^{N}\Rightarrow
X$.)  Consequently, to complete the proof that
\begin{equation}\label{eq:LnObj}
	L_{N} \Longrightarrow L=\exp
 	      \left\{\iint \theta (t,x)\,M(dt,dx) -
	      \frac{1}{2} \iint \theta (t,x)^{2}X (t,x)\,dt \, dx\right\},
\end{equation}
it suffices to prove that 
\begin{equation}\label{eq:qvConvergenceSIR}
	\sum_{t}\sum_{x}\Delta^{N}_{t}(x)^{2}\varrho^{N}_{t} (x)^{2}
	\Longrightarrow 
	 \iint \theta (t,x)^{2}X (t,x)\,dt \, dx.
\end{equation}

\begin{proof}
[Proof of \eqref{eq:qvConvergenceSIR}] This uses only
Theorem~\ref{theorem:BRW} and a variance calculation. Define
\begin{align}\label{eq:AN}
	A^{N}&=  \sum_{t}\sum_{x} 
	\Delta ^{N}_{t} (x)^{2}\varrho^{N}_{t} (x)^{2}
	\quad \text{and}\\
\label{eq:BN}
	B^{N}&= \sum_{t}\sum_{x} 
	        \lambda^{N}_{t} (x)\varrho^{N}_{t} (x)^{2}/N^{2\alpha} ;
\end{align}
we will show that 
\begin{equation}\label{eq:ANBN}
	A^{N}-B^{N} =o_{P} (1)
\end{equation}
Since $\lambda^{N}_{tN^{\alpha}}(xN^{\alpha /2})/N^{\alpha /2}- X^{N}
(t,x)=o_{P} (1)$ as $N \rightarrow \infty$, by
Theorem~\ref{theorem:BRW} (recall \eqref{eq:lambda} that
$\lambda^{N}(t,x)$  is the average of the counts $Y^{N}_{t-1} (x')$
over the neighboring sites $x'=x,x\pm 1$), 
\[
	B^{N}=\iint  X^{N} (t,x) \theta^{N} (t,x)^{2} \,dt \,dx +o_{P} (1).
\]
Theorem~\ref{theorem:BRW} also implies that 
\[
	\iint  X^{N} (t,x) \theta^{N} (t,x)^{2} \,dt \,dx 
	\Longrightarrow \iint  X (t,x) \theta (t,x)^{2} \,dt \,dx .
\]
Therefore, proving \eqref{eq:ANBN} will establish \eqref{eq:qvConvergenceSIR}.

For any constant $C<\infty$, define
\[
	\tau_{C}=\tau^{N}_{C}=
	\max \{t\, : \, \max_{x}\max_{s\leq t}Y^{N}_{s-1}(x)\leq CN^{\alpha /2}\},
\]
and let $A^{N}_{C}$ and $B^{N}_{C}$ be the restrictions of the sums
\eqref{eq:AN}--\eqref{eq:BN} to the range $t\leq \tau_{C}\wedge
CN^{\alpha}$ and $|x|\leq CN^{\alpha /2}$.  Note that $\tau_{C}$ is a
stopping time, and that $\lambda^{N}_{t} (x)\leq CN^{\alpha /2}$ on
the event $t\leq \tau_{C}$. Since the range of summation in
\eqref{eq:AN} and \eqref{eq:BN} is limited by \eqref{eq:GWduration}
and \eqref{eq:maxSpatialBRW},  for any $\varepsilon >0$
there exists $C=C_{\varepsilon}<\infty$, independent of $N$, such that
\[
	A^{N}=A^{N}_{C} \quad \text{and} \quad B^{N}=B^{N}_{C}
\]
except possibly on events of probability $\leq \varepsilon$. Thus, it
suffices to prove \eqref{eq:ANBN} with $A^{N},B^{N}$ replaced by
$A^{N}_{C},B^{N}_{C}$. For this, we use the fact that the offspring
counts $Y^{N}_{t} (x)$ have conditional Poisson distributions (given
the past $\mathcal{H}_{t-1}$) with
means $\lambda^{N}_{t} (x))$: This implies that the conditional means
and conditional variances coincide, and that
\begin{equation}\label{eq:Poisson4th}
	E_{P} ((( Y^{N}_{t} (x)-\lambda^{N}_{t} (x))^{2} -\lambda^{N}_{t} (x))^{2}
	      \, | \, \mathcal{H}^{N}_{t-1})=2\lambda^{N}_{t}(x)^{2}.
\end{equation}
For $t\leq \tau_{C}$, the right side is bounded by
$2C^{2}N^{\alpha}$. Furthermore, for $t\leq \tau_{C}\wedge C$, 
\[
	\varrho^{N}_{t} (x)^{4}\leq C^{8}N^{2\alpha}/N^{4-4\alpha}.
\]
Thus, the conditional variance of each term in the sum
$A^{N}_{C}-B^{N}_{C}$ is bounded by $2C^{10}N^{7\alpha -4}$.
Since the number of nonzero terms in the sum is
$C(2C+1)N^{3\alpha/2}$, it follows that    
\[
	E (A^{N}_{C}-B^{N}_{C})^{2}\leq C' N^{17\alpha /2-4}.
\]
But $\alpha \leq 2/5$, so this converges to $0$ as $N \rightarrow
\infty$. 
\end{proof}
 
\subsection{Likelihood ratios: \emph{SIS} epidemics}\label{ssec:LRsis}
For \emph{SIS} epidemics, Theorem~\ref{theorem:sis} requires that
$\alpha \leq 2/3$, so by \eqref{eq:GWduration} and
\eqref{eq:maxSpatialBRW}, the number of sites/times $(x,t)$ that
contribute nontrivially to the likelihood ratio product
\eqref{eq:likelihoodratio} is $O_{P} (N^{3\alpha /2})=O_{P}
(N)$. Thus, error terms of magnitude $o_{P} (N^{-1})$ can be ignored
in each factor. 

In the modified standard couplings for \emph{SIS} epidemics, the
conditional probability that there is a blue offspring of a red parent
at $(x,t)$, given $y$ red-parent offspring in total at $(x,t)$, is
$\kappa_{N} (y)=y (y-1)/2N$ (see
\eqref{eq:kappaDefinitionSIS}). Hence,
\begin{align*}
	L_{N} (t,x)&=1-\binom{y}{2}/N +\frac{\lambda}{y+1}\binom{y+1}{2}/N\\
	      &=1-y (\lambda -y-1)/2N\\
	      &=1-\lambda (y-\lambda)/2N - (y-\lambda)^{2}/2N +y/2N.
\end{align*}
Since terms of order $o_{P} (N^{-1})$ can be ignored,
\eqref{eq:likelihoodratio} can be written as
\begin{equation}\label{eq:lrSIS}
	L_{N}= (1+o_{P} (1))\prod_{t\geq 1}\prod_{x\in \zz{Z}}\exp \left\{ 
			 -A_{N} (t,x)-B_{N} (t,x)-C_{N} (t,x)/2	   \right\}
\end{equation}
where
\begin{align*}
	A_{N} (t,x)&=\lambda (y-\lambda)/2N;\\
	B_{N} (t,x)&=\lambda^{2} (y-\lambda)^{2}/8N^{2};\\
	C_{N} (t,x)&=(y-\lambda)^{2}/N -y /N.
\end{align*}
Here we continue to use the convention $y=y_{t} (x)$ and $\lambda
=\lambda^{N}_{t} (x)$, as in \eqref{eq:likelihoodratio}. Hence, to
prove the convergence \eqref{eq:LnObj}  $L_{N}\Rightarrow L$ (jointly
with that of $X^{N} (t,x)\Rightarrow X (t,x)$), it will
suffice to show that 
\begin{gather}\label{eq:SISObjA}
	\sum_{t,x} A_{N} (t,x)\Longrightarrow \iint  \theta (t,x)\,M (dt,dx);\\
\label{eq:SISObjB}
	\sum_{t,x}B_{N} (t,x)\Longrightarrow \iint \theta (t,x)^{2}X
	(t,x)\,dt\, dx;
		    \quad \text{and}\\
\label{eq:SISObjC}
	\sum_{t,x} (C_{N} (t,x)
	(t,x))\Longrightarrow  0
\end{gather}
where $\theta (t,x)=0$ for $\alpha <2/3$ and $\theta (t,x)=X (t,x)/2$
for $\alpha =2/3$, as in parts (a)-(b) of Theorem~\ref{theorem:sis}.

\begin{proof}
[Proof of \eqref{eq:SISObjA}]
This  is virtually identical to the proof of
the analogous convergence in the \emph{SIR} case. The increments
$\Delta^{N}_{t} (x):= (Y^{N}_{t} (x)-\lambda^{N}_{t} (x))/N^{\alpha}$
are the  masses \eqref{eq:ommExplicit} in the orthogonal martingale
measures $M^{N}$, and so
\[
	\sum_{t,x} A_{N} (t,x)=\iint \theta^{N} (t,x) \, M^{N} (dt,dx),
\]
where $\theta^{N} (t,x)=\lambda^{N}_{N^{\alpha}t} (xN^{\alpha
/2})/2N^{1-\alpha}$. If $\alpha <2/3$ then $\max_{t,x}\theta^{N} (t,x)
\rightarrow 0$ in probability, whereas if $\alpha =2/3$ then
$\theta^{N} (t,x)$ meets the requirements of
Corollary~\ref{corollary:SIConvergenceSISCase}.Thus, the convergence
\eqref{eq:SISObjA} follows from Corollary
\ref{corollary:SIConvergenceSISCase}.
\end{proof}

\begin{proof}
[Proof of \eqref{eq:SISObjB}] This proceeds in the same manner as the
proof of \eqref{eq:qvConvergenceSIR} in the \emph{SIR} case, by
showing that the term $B_{N} (t,x)$ in the sum \eqref{eq:SISObjB} can
be replaced by $B_{N}' (t,x):=\lambda^{3}/8N^{2}$. To do this, we
truncate the sum of the differences in exactly the same way as in the
proof of \eqref{eq:qvConvergenceSIR}, using the same stopping times
$\tau_{C}=\tau^{N}_{C}$. Note that the truncated sum once again has $C
(2C+1)N^{3\alpha/2}$ terms, and that in each term, $\lambda \leq
CN^{\alpha /2}$. Using the conditional variance
formula~\eqref{eq:Poisson4th}, we find that after truncation,
\[
	E \left(\sum \sum (B_{N} (t,x)-B_{N}' (t,x)) \right)^{2}
	\leq C'N^{3\alpha /2}N^{6\alpha /2}/N^{4}
	\leq C'N^{-1},
\]
since $\alpha \leq 2/3$. Consequently, $\sum \sum B_{N}$ can be
replaced by $\sum \sum B_{N}'$. But
\[
	\sum_{t} \sum_{x}B_{N}' (t,x) =
	\sum_{t} \sum_{x} \lambda^{N}_{t} (x)^{3}/8N^{2}
	\Longrightarrow 
	\iint X (t,x)\theta (t,x)^{2}  \,dx \,dt
\]
by Theorem~\ref{theorem:BRW} and the definition of $\theta$ (parts
(a)--(b) of Theorem~\ref{theorem:sis})
\end{proof}

\begin{proof}
[Proof of \eqref{eq:SISObjC}] 
This is  based on a variance calculation similar to  those used to
prove \eqref{eq:SISObjB} and \eqref{eq:qvConvergenceSIR}. Note first
that the terms $C_{N} (t,x)$ constitute a martingale difference
sequence relative to the  filtration $\mathcal{H}_{t}$, because the
mean and variance of a Poisson random variable coincide. The
conditional variances of the terms $C_{N} (t,x)$ can be estimated
as follows: If $Y\sim$Poisson$(\lambda)$ with $\lambda \geq 1$, then
\begin{align*}
	E ((Y-\lambda)^{2}-Y)^{2}&=
	E ((Y-\lambda)^{2}-\lambda )^{2}+
	E (Y-\lambda)^{2}+2E((Y-\lambda)^{2}-\lambda )(Y-\lambda)\\
	&\leq 2\lambda^{2}+\lambda +2\sqrt{2\lambda^{2}\lambda }\\
	&\leq 5\sqrt{2}\lambda^{2}.
\end{align*}
Now truncate the sum $\sum \sum C_{N} (t,x)$ as before (that is,
$t\leq \tau_{C}\wedge CN^{\alpha}$ and $|x|\leq CN^{\alpha /2}$):
Since the number of terms is $C'N^{3\alpha /2}$ and $\lambda \leq
CN^{\alpha /2}$ for each term, the
variance of the truncated sum is bounded by
$C''N\lambda^{2}/N^{2}$. For $\alpha \leq 2/3$, this converges to zero
as $N \rightarrow \infty$.
\end{proof}

\section{ Weak Convergence in $\zz{D} ([0,\infty ),C_{b} (\zz{R}))$:
Proof of Theorem \ref{theorem:BRW} }\label{sec:tightness}

In this section we prove
Theorem~\ref{theorem:BRW} by verifying that the sequence of random
functions $X^{M} (t,x)$ is tight, provided the hypotheses of
Theorem~\ref{theorem:BRW} on the initial densities hold.

\subsection{Moment Estimates}\label{ssec:momentEstimates}
The proof will be  based on moment estimates for occupation counts of a
branching random walk started by a single particle located at the
origin at time $0$. Denote by $Y_{n} (x)$ the number of particles at
site $x\in \zz{Z}$ at time $n\in \zz{Z_{+}}$, and by $E^{x},P^{x}$ the
expectation operator and probability measure under which the branching
random walk is initiated by a single particle located at $x$. For
notational ease set
\begin{align}\label{eq:gaussianDensity}
	\varphi (x)&=\exp \{-x^{2} \},\\
\notag 
	\varphi_{n} (x)&=\exp \{-x^{2}/n \},\\
\notag 
	\Phi_{n} (x,y)&=\varphi_{n} (x)+\varphi_{n} (y).
\end{align}

\begin{proposition}\label{proposition:momentEstimates}
There exist finite constants $C_{m},\beta_{m}$ 
such that  for all $m,n\in \zz{N}$, all $x,y\in \zz{Z}$, and all $\alpha \in (0,1)$
\begin{align}\label{eq:EstA}
	E^{0}Y_{n} (x)^{m}&\leq C_{m}n^{-1}n^{m/2} 
	\varphi_{n} (\beta_{m} x),\\
\label{eq:EstB}
	|E^{0} ( Y_{n} (x)-Y_{n} (y))^{m}|
	&\leq C_{m}n^{-1}n^{m/2}| (x-y)/\sqrt{n}|^{m/5} \Phi_{n} (\beta x,\beta y)
	     ,\\
\label{eq:EstC} 
		|E^{0} ( Y_{n} (x)-Y_{n+[\alpha n]} (x))^{m}|
		&\leq C_{m}n^{-1}n^{m/2}
		\alpha^{m/5}\varphi_{n} (\beta_{m} x).
\end{align}
\end{proposition}

The (somewhat technical) proofs will be given in sections
\ref{ssec:EstA}--\ref{ssec:EstC} below.  The exponents $m/5$ on
$\alpha$ and $|x-y|/\sqrt{n}$ in \eqref{eq:EstC} and \eqref{eq:EstB}
are not optimal, but the orders of magnitude in the estimates make
sense, as can be seen by the following reasoning: The probability that
a branching random walk initiated by a single particle survives for
$n$ generations is $O (n^{-1})$, and on this event the total number of
particles in the $n$th generation is $O(n)$. Thus, on the event of
survival to generation $n$, the number of particles $Y_{n} (x)$ at a
site $x$ at distance $\sqrt{n}$ from the origin should be $O
(\sqrt{n})$. This is consistent with $m$th moment of order $O
(n^{-1}n^{m/2})$.

\subsection{Tightness in $\zz{D} ([\delta ,\infty ),C_{b}
(\zz{R}))$}\label{ssec:tightness} The proof of tightness will be
broken into two parts: In this section, we will show that for any
$\delta >0$, under the hypotheses of Theorem~\ref{theorem:BRW}, the
density processes $X^{k} (t,x)$ restricted to the time interval
$[\delta,\infty)$ are tight. Using this and an auxiliary smoothness
result of Shiga \cite{shiga:contrasts} for the Dawson-Watanabe
process, we will then conclude in sec.~\ref{ssec:tightness2} that the
density processes $X^{k} (t,x)$ with $t$ restricted to some interval
$[0,\delta]$ are tight.

\bigskip \noindent 
\textbf{Standing Assumptions:} In sections \ref{ssec:tightness}
--\ref{ssec:tightness2}, $Y^{k}_{t} (x)$ will be a sequence of
branching random walks satisfying the hypotheses of Watanabe's
theorem, and $X^{k} (t,x)$ will be the corresponding renormalized
density processes, defined by equation \eqref{eq:Yn}. In addition,
assume that the  initial configurations are such that for every $k$, the
initial particle density function $X^{k} (0,\cdot)$ has support
contained in $J$, for some fixed compact interval $J$.

\begin{proposition}\label{proposition:kolmogorov-chentsov}
There exist constants $C=C_{n},p<\infty$ and $\beta >0$ such that for
every $k=1,2,\dotsc$, all $x,y \in \zz{R}$, and all $s,t \in [1/n,n]$,
\begin{align}\label{eq:objTightnessA}
	E|X^{k} (t,x)-X^{k} (t,y)|^{p}&\leq C|x-y|^{2+\delta } 
	(\varphi (\beta x)+\varphi (\beta y))\quad \text{and}\\
\label{eq:objTightnessB}
	E|X^{k} (t,x)-X^{k} (s,x)|^{p}&\leq C|t-s|^{2+\delta }
	\varphi (\beta x).
\end{align}
\end{proposition}

\noindent \textbf{Note.} Under the conditions of Watanabe's theorem,
the bounds \eqref{eq:objTightnessA}--\eqref{eq:objTightnessB} will in
general hold only for $t,s$ away from zero: in fact, if $X_{0}$ is
singular then the densities $X^{k} (t,x)$ will blow up as $t\rightarrow 0$.

\begin{proof}
The inequalities \eqref{eq:objTightnessA}--\eqref{eq:objTightnessB}
follow from inequalities \eqref{eq:EstB}--\eqref{eq:EstC},
respectively, with $p=m\geq 12$. To see this, observe that the moments in
\eqref{eq:objTightnessA}--\eqref{eq:objTightnessB} can be related to
the corresponding moments for branching random walks $Y$ started from
single particles located at points $x\in J\sqrt{k}$: If
$n_{k}:=\xnorm{Y^{k}_{0}}$ is the number of particles in the initial
configuration of the $k$th \emph{BRW} then for any even integer $m$,
\begin{equation}\label{eq:densityBound}
	E|X^{k} (t,x)-X^{k} (t,y)|^{m}=
	\sum_{i=1}^{m}\sum_{\mathbf{m}\in \mathcal{P}_{r} (m)}
	\sum_{x_{1},x_{2},\dotsc x_{r}}
	\prod_{i=1}^{r}k^{-m_{i}/2}
	 E^{x_{i}} (Y_{[tk]} (\sqrt{k}x)-Y_{[kt]} (\sqrt{k}y))^{m_{i}}
\end{equation}
where $\mathcal{P}_{r} (m)$
is the set of all integer partitions of $m$ with $r$ nonzero elements
$m_{i}$, the inner sum is over all choices of $r$ particles from the
$n_{k}$  particles in the initial
configuration,  and $x_{i}$ is the location of the $i$th
particle. Since $\mathcal{F}_{k}Y^{k}_{0}\Rightarrow X_{0}$, the
masses $n_{k}$ must be asymptotically proportional to  $k$, so the
number of choices in the inner sum is $\leq Ck^{r}$, for some constant
$C$ independent of $k$. Since the initial configurations are all
restricted to lie in $\sqrt{k}J$, with $J$ compact, the bounds implied
by \eqref{eq:EstB} for $E^{x_{i}}$ are comparable to those for
$E^{0}$, after changing $\beta_{m}$ to $\beta_{m}/2$, because the
Gauss kernel $\varphi$ satisfies $\varphi (x-y)\leq C_{J}\varphi (x/2)$
for all $x\in \zz{R}$ and $y\in J$. Thus, by \eqref{eq:EstB},
\[
	RHS\eqref{eq:densityBound}
	\leq C \sum_{i=1}^{m}\sum_{\mathbf{m}\in \mathcal{P}_{r} (m)}
	k^{r} \prod_{i=1}^{r} k^{-1}|x-y|^{m_{i}/5}
	\Phi_{1} (\beta x,\beta y),
\]
provided $t$ is bounded away from $0$ and $\infty$. This clearly
implies \eqref{eq:objTightnessA}. A similar argument gives
\eqref{eq:objTightnessB}. 
\end{proof}

\begin{corollary}\label{corollary:tightness}
Assume that all of the measures $\mathcal{F}_{k}Y^{k}_{0}$ have
supports $\subset J$, where $J$ is compact, and assume that the
hypothesis \eqref{eq:initialCondition} of Watanabe's theorem
holds. Then for every $\delta >0$, the random functions $\{ X^{k}
(t,x)\}_{t\geq \delta}$ converge weakly in the Skorohod space $\zz{D}
([\delta ,\infty ),C_{b} (\zz{R}))$ to the Dawson-Watanabe density
process $\{X (t,x) \}_{t\geq \delta}$ restricted to time $t\geq
\delta$.
\end{corollary}

\begin{remark}\label{remark:2}
Observe that Corollary \ref{corollary:tightness} holds even for
initial conditions $Y^{k}_{0} (\cdot)$  whose Feller-rescalings
$\mathcal{F}_{k}Y^{k}_{0}$ converge to \emph{singular}
measures. 
\end{remark}

\begin{proof}
Since the associated measure-valued processes converge to the
Dawson-Watanabe process, it suffices to show that the sequence
$\{X^{k} (t,x)\}_{t\geq \delta}$ is tight.  This follows from inequalities 
\eqref{eq:objTightnessA}--\eqref{eq:objTightnessB} by the usual
Kolmogorov--Chentsov argument. (See, for instance, \cite{billingsley},
Th.~12.3, or \cite{karatzas-shreve}, Problem~4.11 for the
one-parameter case. Here, since there are two parameters $t,x$, the
exponent must be $>2$.)
\end{proof}

\subsection{Tightness in $\zz{D} ([0,\infty ),
C_{b}(\zz{R}))$}\label{ssec:tightness2} It remains to show that under
the stronger hypotheses of Theorem~\ref{theorem:BRW} the rescaled
particle densities $X^{k}(t,x)$ are tight for $t\in [0,\infty )$. It
is possible to do this by estimating moments, but this is
messy. Instead, we will use a soft argument, based on an estimate for
the Dawson-Watanabe process proved by Shiga (\cite{shiga:contrasts},
Lemma~4.2):

\begin{lemma}\label{lemma:shiga}
Let $X (t,x)$ be the density of a standard Dawson-Watanabe process
with initial condition $X (0,x)=f (x)$ and variance parameter
$\sigma^{2}$. Define
\begin{equation}\label{eq:deviation}
	N (t,x)=X (t,x)-G_{t}f (x)
\end{equation} 
where $G_{t}f$ is the convolution of $f$ with the Gaussian density of
variance $\sigma^{2}t$.  For any compact interval $J$ there exist
constants $C_{1},C_{2}$ such that for all $\varepsilon ,T>0$, if
$f\leq \beta \mathbf{1}_{J}$ then
\begin{equation}\label{eq:shigaIneq}
	P\{|N (t,x)|\geq \varepsilon \beta e^{- (T-t)|x|} \; \text{for
	some} \; t<T/2 \;\text{and}\;x\in \zz{R}\}\leq C_{1}\varepsilon^{-24} 
	\exp \{-C_{2}\varepsilon^{2} T^{-1/4}\}.
\end{equation}
Consequently, for any $\alpha >1$, $\varepsilon >0$, and compact
interval $J$ there exists $T>0$ such that if $f\leq \beta
\mathbf{1}_{J}$ then
\begin{equation}\label{eq:shigaConsequence}
	P\{X (t,x)\geq \alpha \beta \;\text{for some} \;
	     t\leq T \;\text{and}\; x\in \zz{R}\}<\varepsilon .
\end{equation}
\end{lemma}

The strategy now is to use Lemma~\ref{lemma:shiga} to deduce a maximal
inequality for the density $X^{k} (t,x)$ of a branching random walk
over a short time interval $t\in [0,\delta]$. For this, we use the
weak convergence result of Corollary~\ref{corollary:tightness}
together with monotonicity of the \emph{BRW} in the initial
condition. Note that adding particles to the initial configuration of
a \emph{BRW} has the effect of augmenting the original \emph{BRW} by
an independent \emph{BRW} initiated by the set of new particles. Thus,
for any two initial particle configurations $Y^{k,A}_{0}$ and
$Y^{k,B}_{0}$ whose discrepancy satisfies
\[
	|Y^{k,A}_{0} (x)-Y^{k,B}_{0} (x)|\leq 
	\beta k\mathbf{1}_{\sqrt{k}J} (x)
\]
there exist coupled branching random walks $Y^{k,A}_{t}$ and
$Y^{k,B}_{t}$ with initial conditions $Y^{k,A}_{0}$ and $Y^{k,B}_{0}$
whose difference is bounded in absolute value by a branching random
walk $Y^{k,C}_{t}$ with initial condition
\[
	Y^{k,C}_{0} (x)=[\beta k]\mathbf{1}_{\sqrt{k}J}.
\]

\begin{lemma}\label{lemma:maxIneqBRW}
For any $\varepsilon >0$ and compact interval $J$ there exist
$\beta,T>0$ such that if
\begin{equation}\label{eq:InitDelta}
	Y^{k}_{0} \leq \beta  \mathbf{1}_{\sqrt{k}J}
	\quad \text{then} \quad 
	P\{\sup_{t\leq T}\sup_{x\in \zz{R}}X^{k} (t,x)\geq
	\varepsilon\}< \varepsilon .
\end{equation}
\end{lemma}

\begin{proof}
Consider first a sequence of branching random walks with initial
densities $X^{k} (0,\cdot)=2\beta \mathbf{1}_{J^{*}}$, where $J^{*}$
is a compact interval containing $J$ in its interior. By
Corollary~\ref{corollary:tightness}, for any $T>0$ the random
functions $X^{k} (T+t,x)$ (with $t\geq 0$) converge to the density
$X(T+t,x)$ of a standard Dawson-Watanabe process with initial density
$X (0,x)=2\beta \mathbf{1}_{J^{*}} (x)$. (Note: The variance parameter
here is $\sigma^{2}=2/3$.) Hence, by Shiga's Lemma~\ref{lemma:shiga},
if $T>0$ is sufficiently small then
\begin{equation}\label{eq:BRWcompare}
	P\{ \min_{x\in J} X^{k} (T,x)\leq \beta\}<\varepsilon 
	\quad \text{and} \quad 
	P\{\sup_{t\in [T,2T]}\sup_{x\in \zz{R}}
		       X^{k} (t,x)\geq 4\beta\}<\varepsilon .
\end{equation}
It  follows that for each $k$, the particle configuration at time
$T$ is such that the renormalized density $X^{k} (T,\cdot)$ exceeds
$\beta \mathbf{1}_{J}$, except on an event of probability
$<\varepsilon$. 

Now consider branching random walks with initial
densities $\tilde{X}^{k} (0,\cdot)=\beta \mathbf{1}_{J}$. By monotonicity in initial
configurations, the corresponding particle density processes $\tilde{X}^{k} (t,x)$ 
are dominated by the density processes $X^{k} (T+t,x)$ of the
preceding paragraph, except on events of probability
$<\varepsilon$. Consequently, by the second inequality in
\eqref{eq:BRWcompare} and the Markov property,
\[
	P\{\sup_{t\leq T}\sup_{x} \tilde{X}^{k} (t,x)\geq 4\beta\} <2\varepsilon .
\]
\end{proof}

\begin{proof}[Proof of Theorem \ref{theorem:BRW}.]  Assume now that
the initial particle densities $f_{k} (\cdot):= X^{k} (0,\cdot)$
satisfy the hypotheses of Theorem~\ref{theorem:BRW}: in particular,
all have support contained in the compact interval $J$, and
$f_{k}\rightarrow f$ uniformly for some continuous function $f$. By
Corollary~\ref{corollary:tightness}, for every $T>0$ the density
processes $X^{k} (T+t,x)$ converge weakly to $X (T+t,x)$, where
$X(t,x)$ is the density process of the Dawson-Watanabe process with
initial density $f$. By Kesten's Theorem~\cite{kesten:brw}, $\varepsilon
>0$ there is a compact interval $J^{*}\supset J$ such that for every
$T>0$ and $k\in \zz{N}$, with probability at least $1-\varepsilon$,
the function $X^{k}(T,\cdot)$ has support contained in $J^{*}$. By
Shiga's Lemma, for any $\varepsilon >0$, if $T>0$ is sufficiently
small,
\[
	P\{\sup_{x} |X^{k} (T,x)-X^{k} (0,x) |\geq \varepsilon\}<\varepsilon .
\]
Hence, by Lemma~\ref{lemma:maxIneqBRW} and a comparison argument, it
follows that  for any $\varepsilon >0$, if $T>0$ is sufficiently
small,
\[
	P\{\sup_{t\leq T}\sup_{x} |X^{k} (T+t,x)-X^{k} (t,x) |\geq
	\varepsilon\}<\varepsilon . 
\]
Tightness of the sequence $\{X^{k} (t,x) \}_{t\geq 0}$ now follows
from  Corollary~\ref{corollary:tightness}.
\end{proof}

\subsection{Proof of Proposition \ref{proposition:momentEstimates}:
Preliminaries}\label{ssec:prop6prelims}

The proofs of the estimates \eqref{eq:EstA}--\eqref{eq:EstB} are based
on a simple recursive formula for the $m$th moment of a linear
functional $\xbrace{Y_{n},\psi}:=\sum_{x}Y_{n} (x)\psi (x)$ of the
particle density at time $n$. Observe that, for any integer $m\geq 1$,
the $m$th power $\xbrace{Y_{n},\psi}^{m}$ is the sum of all possible
products $\prod_{i=1}^{m}\psi (x_{i})$, where $x_{i}$ is the location
of a particle in the $n$th generation (particles may be repeated). For
any such product, the $r\leq m$ particles involved will have a last
common ancestor (LCE), situated at a site $z$ in the $k$th generation,
for some $k\leq n$. There are two possibilities: either there is just $r=1$
particle in the product, in which case it is its own LCE and
$k=n$, or there are at least two particles, in which case the LCE
belongs to a generation $k<n$. Conditioning on the generation $k$ and
location $z$ of the LCE leads to the following formula (for bounded
functions $\psi$):
\begin{equation}\label{eq:momentFormula}
	E^{x}\xbrace{Y_{n},\psi}^{m}
	=\sum_{z}\zz{P}^{n} (x,z)\psi (z)^{m}
	+\sum_{k=0}^{n-1}\sum_{z}\zz{P}^{k} (x,z)
			\sum_{r=2}^{m}\kappa_{r}
			\sum_{\mathbf{m}\in \mathcal{P}_{r} (m)}
			F_{n-k} (\psi ;\mathbf{m};z)
\end{equation}
where 
\begin{equation}\label{eq:tpm}
	\zz{P} (x,y)=\frac{1}{3}\mathbf{1}\{|x-y|\leq 1 \}
\end{equation}
is the transition probability kernel for nearest neighbor random walk
with holding, $\kappa_{r}$ is the $r$th descending factorial moment of
the offspring distribution (the expected number of ways to choose $r$
particles from the offspring of any particle), $\mathcal{P}_{r} (m)$
is the set of all integer partitions of $m$ with $r$ nonzero elements
$m_{i}$, and
\begin{equation}\label{eq:Fn}
	F_{n} (\psi ;\mathbf{m};z):=
	3^{-r}\sum_{j_{i}=0,1,-1}\prod_{i=1}^{r}
		E^{z+j_{i}}\xbrace{Y_{n},\psi}^{m_{i}}.
\end{equation}
Notice that in each term of the second sum in
\eqref{eq:momentFormula}, $r\geq 2$, reflecting the fact that these
terms correspond to final products with $r\geq 2$ distinct particles,
in which the individual particles are repeated $m_{i}$ times. Since
$\sum_{i=1}^{r}m_{i}=m$ and all $m_{i}\geq 1$, it follows that
$m_{i}<m$: this is what makes formula \eqref{eq:momentFormula}
recursive. The formula \eqref{eq:Fn} accounts for the possibility that
the $r$ offspring will jump to random sites adjacent to the location
$z$ of the LCE.  The appearance of powers $\zz{P^{n}}$ of the
transition probability kernel of the nearest neighbor random walk with
holding derives from the fact that the branching random walk is
\emph{critical}, so that the expected number of descendants at $(n,z)$
of an initial particle at $(0,x)$ is $\zz{P}^{n} (x,z)$. Following are
standard estimates that will be used to bound such transition
probabilities.

\begin{lemma}\label{lemma:LLT}
There exist constants $C,\beta <\infty$ such that for all $x,y\in
\zz{Z}$,  $n\in \zz{N}$, and $0 \leq \alpha \leq 1$,
\begin{align}\label{eq:lltA}
	\zz{P}^{n} (0,x)&\leq Cn^{-1/2}\varphi_{n} (\beta x),\\
\label{eq:lltB}
	|\zz{P}^{n} (0,x)-\zz{P}^{n} (0,y)|&\leq 
		   C ( n^{-1/2}|x-y|\wedge 1)\Phi_{n} (\beta x,\beta y),\\
\label{eq:lltC}
 	|\zz{P}^{n} (0,x)-\zz{P}^{n+\alpha n} (0,x)|&\leq 
		    Cn^{-1/2}\alpha \varphi_{n} (\beta x).
\end{align}
\end{lemma}

\begin{proof}
The first inequality follows from the local limit theorem and standard
large deviations estimates for simple nearest-neighbor random
walk. The second and third inequalities use also the fact that the
Gauss kernel $\varphi (x)$ is uniformly Lipshitz in $x$.
\end{proof}

Finally, we record some elementary inequalities for Gaussian densities:

\begin{lemma}\label{lemma:gaussianConvolution}
For any $\beta>0$ there exists $C<\infty$ such that for all $x\in
\zz{Z}$,  $n\geq 1$,
and $k\geq n/2$,
\begin{gather}\label{eq:gaussianSums}
	\sum_{y\in \zz{Z}}\varphi_{n} (\beta y)
	\leq C\sqrt{n};\\
\label{eq:gaussianConvolution}
	\sum_{y\in \zz{Z}}\varphi_{k} (\beta y)\varphi_{n-k} (\beta
	x-\beta y) \leq C \sqrt{n-k}\varphi_{n} (\beta x/4);\\
\label{eq:betaChange}
	\exp \{-\beta x^{2}/4n \}\leq 
	 C_{\beta}\exp \{-\beta (x \pm 1)^{2}/2n \}.
\end{gather}
\end{lemma}

\subsection{Proof of \eqref{eq:EstA}}\label{ssec:EstA} The proofs of
the inequalities in Proposition~\ref{proposition:momentEstimates} will
proceed by induction on the power $m$. In cases
\eqref{eq:EstA}--\eqref{eq:EstB} , the starting point will be formula
\eqref{eq:momentFormula}; we will use the induction hypothesis to
bound the factors in the products \eqref{eq:Fn}.  In each case it will
be necessary to analyze terms separately in the ranges $k\leq n/2$ and
$k>n/2$. To prevent a proliferation of subscripts, we will adopt the
convention that values of constants $C,\beta$ may change from one line
to the next. In particular, in each inductive step we will relax the
constant $\beta$ (from $\beta $ to $\beta /2$ or $\beta /4$) to
account for differences $\pm 1$ in the arguments of exponentials: this
is justified by \eqref{eq:betaChange} in Lemma
\ref{lemma:gaussianConvolution} above.

For the proof of \eqref{eq:EstA}, use the function $\phi =\delta_{x}$
in formula \eqref{eq:momentFormula}. When $m=1$, the terms indexed by
$0\leq k<n$ in \eqref{eq:momentFormula} all vanish, leaving 
$E^{0}Y_{n} (x)=\zz{P}^{n} (0,x)$. Thus, the inequality
\eqref{eq:EstA} follows for $m=1$ directly  from Lemma~\ref{lemma:LLT}. 
Assume now that  \eqref{eq:EstA} holds for all powers $<m$, where $m\geq 2$.
By \eqref{eq:momentFormula} and the induction hypothesis,
\[
	E^{0}Y_{n} (x)^{m}\leq 
	\zz{P}^{n} (0,x)+C\sum_{k=0}^{n-1}\sum_{z}\zz{P}^{k} (0,z)
		   \sum_{\mathbf{m}\in \mathcal{P} (m)}
		   \prod_{i=1}^{r} ( (n-k)^{-1+m_{i}/2} 
		   \varphi_{n-k}(\beta (x-z)))
\]
where $\mathcal{P} (m)=\cup_{i=2}^{m}\mathcal{P}_{r} (m)$ is the set
of all partitions of $m$ with at least two nonzero elements
$m_{i}$. (Note that the constants $\kappa_{r}$ have been absorbed in
$C$.)  The initial term $\zz{P}^{n} (0,x)$ has already been disposed in
the case $m=1$ (since the right side of the inequality
\eqref{eq:EstA} is nondecreasing in $m$, provided $C_{m}\uparrow$ and
$\beta_{m}\downarrow$).  Thus, we need only consider the terms $k<n$
of the second sum.

Consider first the terms $k\leq n/2$. For indicies in this range,
$n-k>n/2$, and so the factors in the inner products can be handled by
simply replacing each $n-k$ by $n$ (at the cost of a constant
multiplier): Since the number $r$ of factors in each product is at
least $2$, and since the exponents $m_{i}$ in eaach interior product
sum to $m$,
\begin{align*}
	\sum_{k\leq n/2}&\leq Cn^{-2+m/2}\sum_{k\leq n/2}\sum_{z}
		    \zz{P}^{k} (0,z)\varphi_{n} (\beta (x-z))\\
		    &\leq C n^{-1+m/2} \varphi_{n}(\beta x).
\end{align*}

Now consider the terms $k>n/2$. For such terms, $n-k$ is no longer
comparable to $n$, and so the factors in the interior products cannot
be estimated in the same manner as in the case $k\leq n/2$.  However,
when $k>n/2$ the transition probability $\zz{P}^{k} (0,z)$ can be
estimated using the local limit theorem \eqref{eq:lltA}: hence, using
the convolution inequality \eqref{eq:gaussianConvolution} of
Lemma~\ref{lemma:gaussianConvolution} and the induction hypothesis,
and  once again that each interior product has at least two factors,
\begin{align*}
	\sum_{k>n/2}&\leq Cn^{-1/2}\sum_{k>n/2}\sum_{z}
			 \varphi_{k} (\beta z) 
			 (n-k)^{-2+m/2}\varphi_{n-k} (\beta (x-z))\\
			 &\leq Cn^{-1/2}\varphi_{n} (\beta x)
			 \sum_{k>n/2}(n-k)^{-2+m/2}\sqrt{n-k}\\
			 &\leq  Cn^{-1/2}\varphi_{n} (\beta x)
			 n^{m/2-1/2},
\end{align*}
as desired. This proves inequality \eqref{eq:EstA}.
\qed 

\subsection{Proof of \eqref{eq:EstB}}\label{ssec:EstB}
This is also by induction on $m$. For notational ease, set
$\psi_{xy}=\delta_{x}-\delta_{y}$, where $\delta_{z}$  is the
Kronecker delta, and set
\begin{equation*}
	 d_{n} (x,y)= (|x-y|/\sqrt{n})\wedge 1.
\end{equation*}
Consider first the case $m=1$. By formula \eqref{eq:momentFormula} and
the local limit bound  \eqref{eq:lltB},
\begin{equation}\label{eq:Case1}
	|E^{0}\xbrace{Y_{n},\psi_{xy}}|
	=|\zz{P}^{n} (0,x)-\zz{P}^{n} (0,y)|
	\leq Cn^{-1/2}d_{n} (x,y)\Phi_{n} (\beta x,\beta y)
\end{equation}
for suitable constants $C,\beta$.  Inequality
\eqref{eq:EstB} for $m=1$ follows easily.

Assume then that inequality \eqref{eq:EstB} is valid for all positive
integer exponents less than $m$.  The first sum on
RHS\eqref{eq:momentFormula} (the terms with $k=n$) can be bounded
above by $Cn^{-1/2} (\varphi_{n} (\beta x)+\varphi_{n} (\beta y))$
using Lemma~\ref{lemma:LLT}, and this in turn is bounded above by
RHS\eqref{eq:EstB}.  Thus, we need only consider the second sum on
RHS\eqref{eq:momentFormula} (the terms $k<n$). In each of these terms,
the interior products \eqref{eq:Fn} have at least two factors, each
with $m_{i}\geq 1$, and in each of these products the sum of the
exponents $m_{i}$ is $m$. Consequently, by \eqref{eq:momentFormula}
and the induction hypothesis,
\[
	\sum_{k=0}^{n-1}
	\leq C \sum_{k=0}^{n-1}\sum_{z}\zz{P}^{k} (0,z)
	|x-y|^{m/5} (n-k)^{-2+m/2-m/10}\Phi_{n-k} (\beta x-\beta z,\beta y-\beta z).
\]
Consider first the terms in the range $k\leq n/2$: for these terms,
$n-k\geq n/2$, and so $\Phi_{n-k}$ is comparable (after a relaxation
of $\beta $) to $\Phi_{n}$. Hence,
\begin{align*}
	\sum_{k=1}^{n/2} &\leq Cn^{-2}n^{m/2-m/10}|x-y|^{m/5}
		    \sum_{k=1}^{n/2}\sum_{z}\zz{P}^{k} (0,z)
		    \Phi_{n-k} (\beta x-\beta z,\beta y-\beta z)\\
		    &\leq Cn^{-1}n^{m/2-m/10}|x-y|^{m/5}
		    \Phi_{n} (\beta x,\beta y),
\end{align*}
as desired. Now consider the range $k>n/2$: Here we use the local
limit estimate for $\zz{P}^{k} (0,z)$ and the Gaussian convolution
inequality \eqref{eq:gaussianConvolution}  to obtain
\begin{align*}
	\sum_{k=n/2}^{n-1}&\leq C|x-y|^{m/5} \sum_{k=n/2}^{n-1}\sum_{z}
			  \zz{P}^{k} (0,z) \Phi_{n-k} (\beta x-\beta
			  z,\beta y-\beta z)
			  (n-k)^{-2+m/2-m/10}\\
			  &\leq Cn^{-1/2}|x-y|^{m/5} \Phi_{n} (\beta x,\beta y)
			  \sum_{k=n/2}^{n-1}(n-k)^{-3/2+m/2-m/10}\\
			  &\leq Cn^{-1}n^{m/2-m/10}|x-y|^{m/5} \Phi_{n} (\beta x,\beta y).
\end{align*}
The final inequality relies on the fact that the exponent
$-1/2+m/2-m/10$ is positive for all $m\geq 2$. This
proves \eqref{eq:EstB}.
\qed 

\subsection{Proof of \eqref{eq:EstC}}\label{ssec:EstC}
The moment formula \eqref{eq:momentFormula} can no longer be used,
since the expectation on LHS\eqref{eq:EstC} involves the state of the
branching random walk at two different times. However, it is not
difficult to derive an analogous formula: For any integer $m\geq 1$,
the  $m$th  power $(Y_{t} (x)-Y_{s} (x))^{m}$ is a sum of products
$\prod_{i=1}^{m}\psi (x_{i})$, where each $x_{i}$ is the location of a
particle in either the $t$th or $s$th generation, and $\psi
(x_{i})=\pm \delta_{x} (x_{i})$, with the sign $\pm$ depending on
whether the particle is in the $t$th or $s$th generation. As in
formula \eqref{eq:momentFormula}, the particles involved in any such
product must have a last common ancestor  in some generation $k$
before the $s\wedge t$th generation. Conditioning on the generation
and location of this last common ancestor leads to the formula
\begin{align}\label{eq:newMomentFormula}
	E^{y}(Y_{s} (x)-Y_{t} (x))^{m}
	=& \zz{P}^{s} (y,x)+ (-1)^{m}\zz{P}^{t} (y,x)\\
\notag 
	&+\sum_{k=0}^{s\wedge t}\sum_{z}\zz{P}^{k} (y,z)\sum_{r=2}^{m}
	\kappa_{r}\sum_{\mathbf{m}\in \mathcal{P}_{r} (m)}
	G_{s-k,t -k} (x-z;\mathbf{m})
\end{align}
where $\kappa_{r}$ and $\mathcal{P}_{r} (m)$ have the same meanings as
in the moment formula \eqref{eq:momentFormula} and
\begin{equation}\label{eq:gst}
		G_{s,t} (z;\mathbf{m})
		=3^{-r}\sum_{j_{i}=0,1,-1}\prod_{i=1}^{r}
		E^{j_{i}} (Y_{s} (z)-Y_{t} (z))^{m_{i}}.
\end{equation}

We will use \eqref{eq:newMomentFormula} to prove \eqref{eq:EstC} by
induction on the power $m$, using arguments similar to those used in
proving \eqref{eq:EstA} and \eqref{eq:EstB}.  Consider first the case
$m=1$: in this case the last sum in \eqref{eq:newMomentFormula}
vanishes, leaving
\begin{equation}\label{eq:m1EstC}
	E (Y_{n} (x)-Y_{n+\alpha n} (x))=
	\zz{P}^{n} (0,x)-\zz{P}^{n+\alpha n} (0,x).
\end{equation}
Hence, inequality \eqref{eq:EstC} follows immediately from estimate
\eqref{eq:lltC}. 

Assume now that \eqref{eq:EstC} holds for all positive integer
exponents smaller than $m$, for some integer $m\geq 2$.  To prove that
\eqref{eq:EstC} holds for the exponent $m$, consider
RHS\eqref{eq:newMomentFormula}, with $s=n$ and $t=n+\alpha n$.  The
sum of the first two terms coincides with \eqref{eq:m1EstC}. To verify
that this sum, in absolute value, is smaller than RHS\eqref{eq:EstC},
observe that by \eqref{eq:lltC} the sum is smaller than $C\alpha
n^{-1/2}\varphi_{n} (\beta x)$; since $\alpha n$ must be an integer,
$\alpha \geq n^{-1}$, and so
\[
	n^{1/2}\alpha \leq n^{m/2}\alpha^{m/8}
\]
for all $m\geq 2$, as required by \eqref{eq:EstC}. Thus, it remains
to prove that the sum of the terms with $r\geq 2$ in
\eqref{eq:newMomentFormula} is also bounded by RHS\eqref{eq:EstC}.

Note first that it suffices to consider values of $\alpha \leq 1/2$,
because for $1/2\leq \alpha \leq 1$ the factor $\alpha^{m/5}$ on
RHS\eqref{eq:EstC} is bounded below. Now split the terms of the last
sum in \eqref{eq:newMomentFormula} into three ranges: first $k\leq
n/2$, then $n/2\leq k\leq n-n\alpha$, and finally $n-n\alpha \leq k<
n$. In the range $k\leq n-n\alpha$, the induction hypothesis applies,
because for these terms $\alpha '=n\alpha / (n-k)\leq 1$. (Recall that
one of the hypotheses of Proposition~\ref{proposition:momentEstimates}
is that $\alpha \leq 1$, so to use \eqref{eq:EstC} in an induction
argument the implied $\alpha'$ cannot exceed $1$). For $k\leq n/2$,
the ratio $n/ (n-k)$ is bounded above by $2$; consequently,
\begin{align*}
	\sum_{k\leq n/2}&\leq C
	\sum_{k\leq n/2}\sum_{z}
	\zz{P}^{k}(0,z) n^{-2+m/2}\alpha^{m/5}
	\varphi_{n-k} (\beta x-\beta z)\\
	&\leq 
	 Cn^{-2+m/2}\alpha^{m/5}
	 \sum_{k\leq n/2}\varphi_{n} (\beta x)\\
	 &\leq Cn^{-1+m/2}\alpha^{m/5}\varphi_{n} (\beta x),
\end{align*}
which agrees with RHS\eqref{eq:EstC}. Next, consider terms in the
range $n/2\leq k\leq n-n\alpha$: By the induction hypothesis,
the local limit bound \eqref{eq:lltA},
and the Gaussian convolution inequality \eqref{eq:gaussianConvolution},
\begin{align*}
	\sum_{k=n/2}^{n-n\alpha}
	&\leq C\sum_{k=n/2}^{n-n\alpha}\sum_{z}
	\zz{P}^{k} (0,z) (n-k)^{-2+m/2}
		     (\alpha n/ (n-k))^{m/5}\varphi_{n-k} (\beta (x-z))\\
	&\leq C n^{-1/2}n^{m/5}\alpha^{m/5}\varphi_{n} (\beta x)
	\sum_{k=n/2}^{n-n\alpha} (n-k)^{-3/2+m/2-m/5}\\
	 &\leq C n^{-1/2}n^{m/5}\alpha^{m/5}\varphi_{n} (\beta x)
	\sum_{j=1}^{n/2}j^{-3/2+m/2-m/5}\\
	&\leq Cn^{-1+m/2}\alpha^{m/5}\varphi_{n} (\beta x).
\end{align*}
The last inequality uses the fact that the exponent $-1/2+m/2-m/5$ is
positive for all $m\geq 2$. 

Finally, consider the terms in the range $n-n\alpha <k\leq n$. Here
the induction hypothesis cannot be used, because $n\alpha /n-k>1$.
Instead we use the bound \eqref{eq:EstA}, which we have already proved
is valid for all $m$. This, together with \eqref{eq:lltA} and
\eqref{eq:gaussianConvolution},   implies
\begin{align*}
	\sum_{k=n-n\alpha}^{n-1}
	&\leq C \sum_{k=n-n\alpha}^{n-1}\sum_{z}
	\zz{P}^{k} (0,z) (n-k)^{-2+m/2}\varphi_{n\alpha } (\beta x-\beta z)\\
	&\leq Cn^{-1/2} \varphi_{n} (\beta x)
	\sum_{k=n-n\alpha}^{n-1}
	(n-k)^{-3/2+m/2}\\
	&\leq Cn^{-1+m/2}\alpha^{-1/2+m/2}\varphi_{n} (\beta x)\\
	&\leq Cn^{-1+m/2}\alpha^{m/5}\varphi_{n} (\beta x).
\end{align*}
This completes the proof of \eqref{eq:EstC}.
\qed

\bigskip \noindent 
\textbf{Acknowledgments.} The author thanks Regina Dolgoarshinnykh and
Xinghua Zheng for useful discussions. 

\bibliographystyle{abbrv} \bibliography{spatial}

\end{document}